\documentclass{article}
\usepackage[utf8]{inputenc}

\usepackage{geometry} 
\usepackage{amssymb}
\usepackage{amsmath}
\usepackage{amsthm}
\usepackage{enumerate}
\usepackage{graphicx, wrapfig}
\usepackage{fullpage}
\usepackage{caption}
\usepackage{subcaption}
\usepackage{lipsum}
\usepackage{algorithm}
\usepackage{algorithmic}
\usepackage[hidelinks]{hyperref}
\usepackage{xcolor}
\usepackage{bbm}
\usepackage{comment}
\usepackage{authblk}

\theoremstyle{definition}

\theoremstyle{definition}

\newcommand{\R}{\mathbb{R}}

\newcommand{\beginsupplement}{%
        \setcounter{table}{0}
        \renewcommand{\thetable}{S\arabic{table}}%
        \setcounter{figure}{0}
        \renewcommand{\thefigure}{S\arabic{figure}}%
     }

\title{An integrative phenotype-structured partial differential equation model for the population dynamics of epithelial-mesenchymal transition}

\author[1,4]{Jules Guilberteau}
\author[2,4]{Paras Jain}
\author[2,5]{Mohit Kumar Jolly}
\author[3,5]{Camille Pouchol}
\author[1,5]{Nastassia Pouradier Duteil}

\affil[1]{Sorbonne Université, CNRS, Université Paris Cité, Inria, Laboratoire Jacques-Louis Lions (LJLL), F-75005 Paris, France.}
\affil[2]{Department of Bioengineering, Indian Institute of Science, Bangalore 560012, India}
\affil[3]{Universit\'e Paris Cité, FP2M, CNRS FR 2036, MAP5 UMR 8145, F-75006 Paris,
France.}
\affil[4]{\tt\small Equally contributed to the work}
\affil[5]{\tt\small To whom correspondence may be addressed: M.K.J.: \href{mailto:mkjolly@iisc.ac.in}{mkjolly@iisc.ac.in}, C.P.: \href{mailto:camille.pouchol@u-paris.fr}{camille.pouchol@u-paris.fr}, N.P.D.: \href{mailto:nastassia.pouradier\_duteil@sorbonne-universite.fr}{nastassia.pouradier\_duteil@sorbonne-universite.fr}}

\begin{document}

\maketitle

\begin{abstract}
    Phenotypic heterogeneity along the epithelial-mesenchymal (E-M) axis contributes to cancer metastasis and drug resistance. Recent experimental efforts have collated detailed time-course data on the emergence and dynamics of E-M heterogeneity in a cell population. However, it remains unclear how different possible processes interplay in shaping the dynamics of E-M heterogeneity: a) intracellular regulatory interaction among biomolecules, b) cell division and death, and c) stochastic cell-state transition (biochemical reaction noise and asymmetric cell division). Here, we propose a Cell Population Balance (Partial Differential Equation (PDE)) based model that captures the dynamics of cell population density along the E-M phenotypic axis due to abovementioned multi-scale cellular processes. We demonstrate how population distribution resulting from intracellular regulatory networks driving cell-state transition gets impacted by stochastic fluctuations in E-M regulatory biomolecules, differences in growth rates among cell subpopulations, and initial population distribution. Further, we reveal that
    %proportionality between heterogeneity and cell growth rates 
    a linear dependence of the cell growth rate on the population heterogeneity
    is sufficient to recapitulate the faster \textit{in vivo} growth of orthotopic injected heterogeneous E-M subclones reported before experimentally. Overall, our model contributes to the combined understanding of intracellular and cell-population levels dynamics in the emergence of E-M heterogeneity in a cell population. 
\end{abstract}

\section{Introduction}
 Intra-tumour heterogeneity – the co-existence of multiple distinct cellular phenotypes in a tumour – is being increasingly reported to associate with poor patient outcomes. It contributes to both metastasis and therapy resistance – two major unsolved clinical challenges \cite{jacquemin2022Dynamic, marusyk2020Intratumor}. Such heterogeneity can arise at a genetic level over the course of tumour evolution and manifests as different clonal populations. However, over the last two decades, non-genetic or phenotypic heterogeneity among cancer cells has been identified as a key driver of disease aggressiveness \cite{bell2020Principles, pillai2023Unraveling}. Such heterogeneity is often characterised by single-cell measurements (flow cytometry, mass cytometry, RNA-seq, ChIP-seq, ATAC-seq),
 showing diversity among cells in a population at the proteome, transcriptome and epigenome 
 levels. A canonical example of non-genetic heterogeneity is along the Epithelial-Mesenchymal (E-M) phenotypic spectrum. Given the implications of E-M heterogeneity in cancer metastasis and patient outcomes, various in vitro, in vivo and in silico attempts have been focused on understanding its underlying mechanisms \cite{brown2022Phenotypic, jain2022Population, karacosta2019Mapping, sahoo2021mechanistic}.
 
An iterative cross-talk among \textit{in silico} and \textit{in vitro} studies has contributed enormously to understanding how non-genetic heterogeneity emerges in a population. Mathematical modelling of regulatory networks underlying Epithelial - Mesenchymal Plasticity (EMP) have reported the co-existence of multiple cellular phenotypes – Epithelial (E), Mesenchymal (M) and one or more hybrid (E/M) cell states \cite{font-clos2018Topography, hari2020Identifying, hong2015Ovol2, rashid2022Network, steinway2015Combinatorial}. Their co-existence has been experimentally reported in varying ratios in multiple cell lines and primary tumours \cite{george2017Survival, karacosta2019Mapping}. 
Synthetic perturbation of the underlying regulatory network led to the loss of bimodal nature of canonical epithelial markers such as E-cadherin and altered the phenotypic distribution \cite{ruscetti2016HDAC, brown2022Phenotypic, celia2018hysteresis, subbalakshmi2020nfatc}. 
The relative stability of cells in different phenotypes, and consequently the phenotypic distribution at a given time, is governed by the underlying topology of the regulatory network involving transcriptional and translational control \cite{hari2022Landscape}, as well as by epigenetic regulation (chromatin modification). Happening at timescales slower than transcriptional regulation, the epigenetic regulation can be reversible or irreversible, as observed experimentally. Thus, epigenetic remodelling can lock cells transiently or permanently in a cell-state, impacting the population phenotypic distribution. Further, cell-to-cell communications either through neighbourhood interactions or paracrine signalling also reshapes the phenotypic distribution \cite{boareto2016Notch, jolly2017Inflammatory, neelakantan2017EMT, yamamoto2017Intratumoral}. 
Hence, a diversity of regulatory interactions within and among cells can contribute to shaping the E-M heterogeneity patterns in a cell population.   

Another milieu of factors such as asymmetric cell division \cite{hitomi2021aAsymmetric, tripathi2020mechanism}, stochastic biochemical noise \cite{munsky2015Integrating}, differences in cellular microenvironment \cite{pally2022Extracellular}, and variable cell-cycle dynamics \cite{lovisa2015Epithelial} can amplify dynamic heterogeneity in a cellular population. These factors alter cell-to-cell variability in protein levels in a population, thus contributing to their functional heterogeneity 
(differential activation of signalling pathways) 
when these cells are exposed to cytotoxic or EMT-inducing growth factors \cite{spencer2009non, strasen2018cell} 
Despite extensive efforts in investigating above-mentioned regulatory and stochastic processes in E-M  plasticity and heterogeneity, only a few computational models have incorporated these processes  within a growing and dividing heterogeneous cellular population. Broadly speaking, two modelling approaches are employed:
\begin{enumerate}
    \item Agent-Based Cell Population models: population models that capture regulatory and stochastic dynamics for individual cells and then generate an ensemble to generate population distribution \cite{jain2022Population, jain2023epigenetic, tripathi2020mechanism}.
    \item Cell Population Balance models: population models that capture the regulatory dynamics and stochastic dynamics %(via statistical distributions)
    of groups of cells with similar cell-states without dealing with individual cell level information. %PDE models, which describe densities of cells. %are  that capture the regulatory dynamics and stochastic dynamics (via statistical distributions) of groups of cells with similar cell-states without dealing with individual cell level information.
\end{enumerate}
%The former approach is more precise since one can follow the fate of each cell, but simulating agent-based models leads to prohibitive computational costs for large number of cells. 
The latter approach has been adopted widely because its output (cell density) can be related directly to flow cytometry experiments conducted at multiple timepoints for a population. These models describe the evolution of the cell density in the phenotypic space. More precisely, they track the number of cells that have a given phenotype in the E-M landscape, where the phenotype is defined by a vector containing the concentration of relevant molecules determining the phenotype of a cell. Further, Cell Population Balance models have been used to combine complex regulatory phenomena, e.g. positive feedback loops \cite{mantzaris2007single}, caspase activation cascade during programmed cell death \cite{hasenauer2011identification}, cell-to-cell communication \cite{shu2011Bistability}, and two mutually inhibiting nodes \cite{spetsieris2009novel}, with stochastic processes like asymmetric cell division and stochastic biochemical reactions for phenotype-structured, and age-structured cell populations \cite{hasenauer2012analysis, schittler2013new, shu2011Bistability}. Finally, since Agent-Based Cell Population model simulate dynamics of each cell individually, they become computationally intractable for realistically large numbers of cells, contrary to models for cell densities. 

Therefore, given its wide applicability, we have here used Cell Population Balance modelling to study the population dynamics of E-M heterogeneity.
Our cell population balance model combines three main cellular processes: 1) \textit{growth} due to cell division and cell death, 2) \textit{state regulation}, which corresponds to the time-evolution of the aforementioned molecules, based on corresponding ODE models and 3) \textit{stochastic cell transition}, which aggregates all sources of stochasticity in the fate of a cell's phenotype.

With the developed model, we demonstrate how heterogeneity along the Epithelial-Mesenchymal axis emerges at a population level. First, we show that the Cell Population Balance model  can capture the previously reported dynamical features of hysteresis, and epigenetic regulation during cells undergoing an Epithelial-Mesenchymal Transition (EMT) followed by a Mesenchymal-Epithelial Transition (MET). Second, we report how cellular heterogeneity depends on the characteristic of stochastic biochemical noise in – 1) external EMT inducing signal (SNAIL), and 2) E-M state variable (miR200 or ZEB levels) – along with the differences in the relative growth rate among E, E/M and M states, and the initial distribution of population in the EMT cell state spectrum. Finally, we show that a population with its fitness (growth rate) proportional to the heterogeneity of population (Renyi entropy) can explain faster tumour growth in vivo and higher proliferation rates in vitro as observed for parental and intermediate clones derived from SUM149T and PMC42-LA cell lines, respectively \cite{bhatia2019Interrogation, brown2022Phenotypic}. 

For all simulations, we use a scheme from the family of particle methods, based on solving a properly defined set of ordinary differential equations (ODEs)~\cite{chertock2017practical}. These are particularly adapted for models based on PDEs such as those developed in the present paper. For the derivation and analysis of particle methods in that context, we refer to the recent theoretical work~\cite{alvarez2023particle}, where it is proved that solutions are properly approximated by the proposed numerical method.

\section{Results}
%\section{Results}

\subsection{Building a Cell Population Balance Model}

Cell Population Balance models allow to compute the evolution of the number of cells of phenotype $y\in \R^n$ at time $t$, denoted by $u(t,y)$.
Here, $n$ represents the dimension of the phenotypic space, that we will also refer to as ``state space''. In our context, a cell's phenotype is considered to be determined by the concentration of several different (internal or external) molecules: the dimension $n$ of the phenotypic space will then correspond to the number of different molecules that determine its phenotype.

%As a first approach, we will consider the simplest case, in which a cell's phenotype is assumed to be determined by the concentration of one given molecule, so that the phenotypic space is one-dimensional ($d=1$). 

To give an intuition for the final PDE model, we consider the simplest case, in which a cell's phenotype is assumed to be determined by the concentration of one given molecule, so that the phenotypic space is one-dimensional ($n=1$). 

The function $u$ is to be understood in the following way: 
for any phenotype interval $[a,b]$, the integral $\int_a^b u(t,y) \,dy$ represents the number of cells whose phenotypes belong to $[a,b]$. Thus, the total number of cells is computed by 
$\int_\R u(t,y) \,dy$, which will be denoted by $\rho(t)$. In our model, the total number of cells $\rho(t)$ is amenable to evolve in time, but for synthetic purposes, from here onward, we will refer to $u$ as the \textit{cell density}.

%We consider the evolution of the number of cells whose phenotype belongs to a small interval $[a,b]$ contained in $\R$ at time $t$, which is computed by the integral $\int_a^b u(t,y) \,dy$. 

In practice, as will be seen in the subsequent sections, 
two molecules, that are markers of EMT, will be taken into account for the description of the (two-dimensional) phenotypic space: miR200 and ZEB. The phenotypic space will then be reduced to one dimension by empirically introducing an artificial state variable~$x$, roughly equivalent to miR200. Then, the cell density $u$ will be expressed in units of cells per number of molecules. As an example, to compute the number of cells which contain fewer than 500 molecules $x$, we compute the integral $\int_0^{500} u(t,y) dy$. With the previous notations, $a=0$ and $b=500$.

\begin{figure}[h!]
    \centering
    \includegraphics[width=\textwidth]{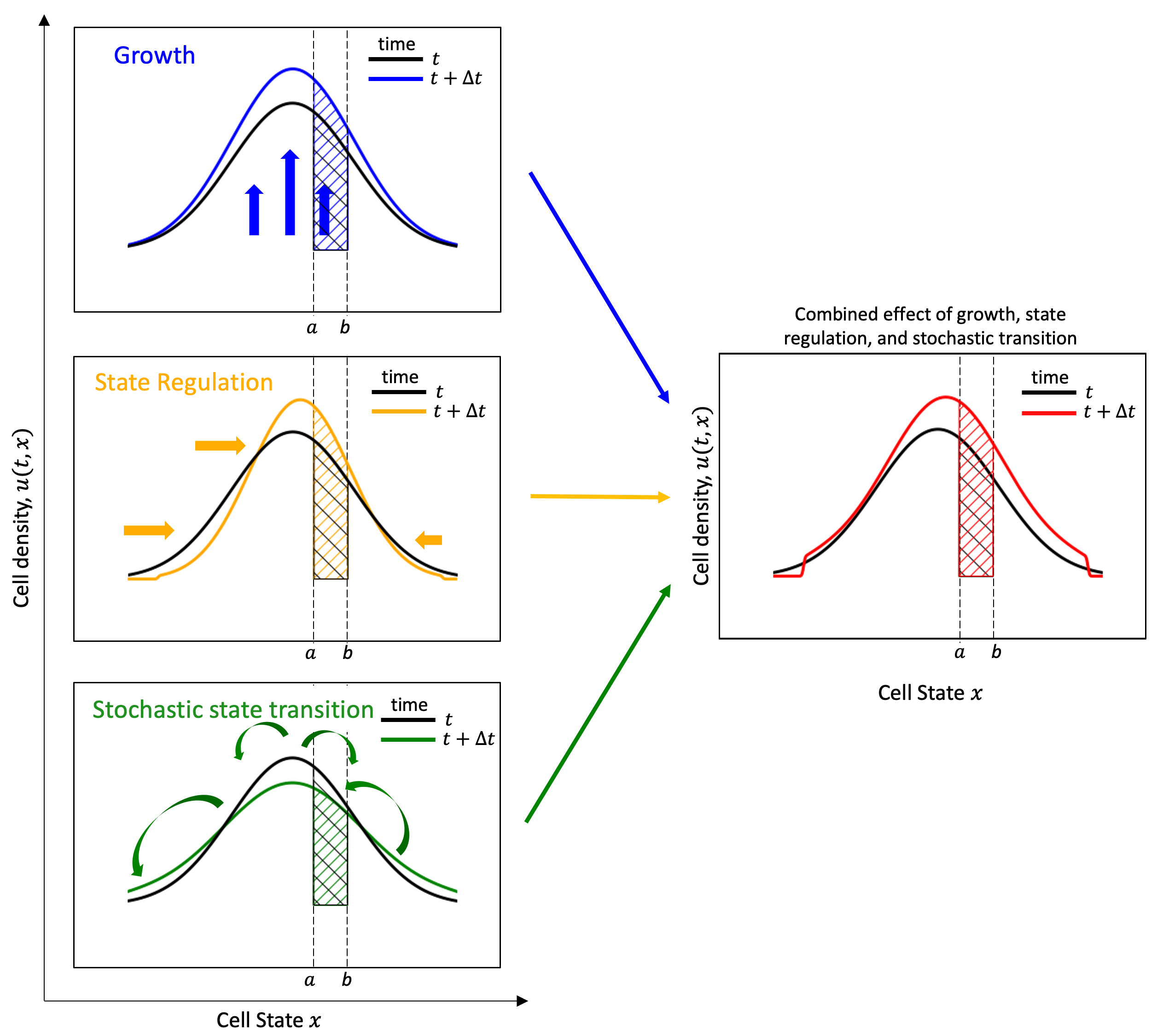}
    \caption{ Growth, state regulation and stochastic state transition combine to give the cell density evolution in the proposed Cell Population Balance model. In each panel, the cell density $u$ is represented at time $t$ in black, and at a subsequent time $t+\Delta t$ in color. The number of cells with cell state in the interval $[a,b]$ is represented by the hashed (black or colored) area. Growth (which includes cell division and cell death) results in an upward or downward shift of the cell density. In this toy example (top left panel), cell division takes place faster than cell death, so the cell density at time $t+\Delta t$ is shifted upward compared to the cell density  at time $t$. Consequently, the number of cells with cell state in $[a,b]$ increases between $t$ and $t+\Delta t$.
    State regulation results in a horizontal displacement of the cell density. In this example (center left panel), a concentration phenomenon is taking place towards the center of the state space. Here, the number of cells with cell state in $[a,b]$ increases between $t$ and $t+\Delta t$.
    Stochastic cell transition is a mixing of the cell population. The evolution of the number of cells with cell state in $[a,b]$ depends on the cell density at every other point of the state space. This results in a diversification of the population, that is in a flattening of the cell density curve (bottom left panel). In this example, the number of cells with cell state in $[a,b]$ decreases between $t$ and $t+\Delta t$.  
    All three mechanisms combine in the cell population balance model to give the evolution of the cell density (right panel).}
    \label{figure 0}
\end{figure}

Three mechanisms will be considered to participate in increasing or decreasing the number of cells within each interval $[a,b]$: growth, cell regulation and stochastic state transition. Their interplay is illustrated in Fig. \ref{figure 0}.

\paragraph{Growth.} The first mechanism that we include in our model is \textbf{growth}, which takes into account
%accounts for new cells of phenotype $y \in \E$ created by 
cell division and cell death. Each cell of phenotype $y$ divides at rate $r(y)$ and its daughter cells are assumed to be given the same phenotype $y$. Hence, the quantity of new cells at time $t+\Delta t$ with phenotype in $[a,b]$ is given by $ \Delta t \int_a^b r(y) u(t,y) \, dy$.
Cells are also assumed to die at a rate $d(y)\rho(t)$, proportional to the total number of cells $\rho(t) = \int_{\R} u(t,y) \,dy$. The quantity of cells that died between the times $t$ and $t+\Delta t$ is then approximated by $ \Delta t \int_a^b d(y) \rho(t) u(t,y) \, dy$. 
The term `$r(y)-d(y)\rho(t)$' thus represents the net growth rate, where $r(y)$ is the intrinsic growth rate, and $d(y)\rho(t)$ the death rate. It is positive if cell division is faster than cell death, and negative in the opposite case.
Overall, the evolution of the number of cells between times $t$ and $t+\Delta t$ caused by the cell population growth is given by 
\begin{align*} \int_a^b u(t+\Delta t,y) \,dy - \int_a^b u(t,y) \,dy   &= \Delta t \int_a^b \left(r(y)-d(y)\rho(t) \right) u(t,y) \, dy.
\end{align*}
Depending on the sign of the right-hand side, the growth mechanism will result in an upward or downward vertical shift of the curve representing the cell density  (see Fig. \ref{figure 0}, top left panel).

\paragraph{State Regulation.} We take advantage of ODE models built to describe the evolution of the phenotype $y$ of one given cell, which takes the form $\dot y = f(y)$.
In the PDE framework, the so-called \textbf{advection} term accounts for all cells whose phenotypes will enter or leave the phenotypic region $[a,b]$ during a small time interval as a result of their inner evolution. Taking only this mechanism into account, the variation of the number of cells in the region $[a,b]$ between two time instants $t$ and $t+\Delta t$ is computed as:
\begin{align*} \int_a^{b} u(t+\Delta t,y) \,dy - \int_a^{b} u(t,y) \,dy   &= \Delta t \; f(a)  u(t,a)- \Delta t \; f(b) u(t,b).
\end{align*}
%Since $\Delta t$ is positive, and $u(t,\cdot)$ is positive (representing a quantity of cells), 
The sign of $f$ at points $a$ and $b$ determines whether cells enter or leave the region. For instance, if $f(a)>0$, the first term of the right-hand side is positive, which translates the fact that cells enter the region $[a,b]$ at point $a$. Similarly, if $f(a)<0$, cells leave the region $[a,b]$ at point $a$.
On the other hand, if $f(b)>0$, cells will leave the region $[a,b]$ at point $b$, and if $f(b)<0$, cells will enter the region $[a,b]$ at point $b$. Intuitively, cells will then be transported towards the right when $f$ is positive, and towards the left when $f$ is negative. Concentration phenomena will happen at points where $f$ is zero and has negative derivative: these points are asymptotically stable points for the ODE $\dot y = f(y)$.
Figure \ref{figure 0} illustrates a toy situation in which $f$ has a stable equilibrium point towards the middle of the phenotypic space, resulting in concentration of the cell density around this point (Fig. \ref{figure 0} left middle panel).

\paragraph{Stochastic cell transition.} The third mechanism that we take into account is stochastic cell-state transition, that we will also refer to as \textbf{mutations}. Here, $M(y,z)$ represents the (infinitesimal) probability that a cell's phenotype $z$ changes to a phenotype $y$. Thus the number of new cells with phenotypes in the interval $[a,b]$ at time $t+\Delta t$ is computed as $\Delta t \int_a^b \Big(\int_{\R} M(y,z) u(t,z)\, dz\Big)\,dy$.
Symmetrically, the number of cells whose phenotypes were in $[a,b]$ at time $t$ and who mutated between times $t$ and $t+\Delta t$ is computed by 
$\Delta t\int_a^b \Big(\int_{\R} M(z,y)\, dz \Big) \, u(t,y) \, dy$. This mechanism is fundamentally non-local, in that the evolution of the number of cells with phenotypes in the interval $[a,b]$ depends on the whole cell population. Often, this mechanism results in a mixing of the population, that is in a flatter cell density (as illustrated in Fig. \ref{figure 0}, bottom left panel). 

\paragraph{Resulting PDE model.}
Putting everything together, the variation of cells whose phenotype belongs to $[a,b]$ during a time interval $\Delta t$ is approximated by 
\begin{align*} \int_a^b u(t+\Delta t,y) \,dy - & \int_a^b u(t,y) \,dy   = \Delta t \Bigg[\underbrace{ \int_a^b \left(r(y)-d(y)\rho(t) \right) u(t,y) \, dy}_{\text{growth}} \\
&+\underbrace{ f(a)  u(t,a)-  f(b) u(t,b)}_{\text{advection}}  \\
&+\underbrace{ \int_a^b \Big(\int_{\R} M(y,z) u(t,z)\, dz\Big)\,dy- \int_a^b \Big(\int_{\R} M(z,y)\, dz \Big) \, u(t,y) \, dy}_{\text{mutations}}  \Bigg].
%+ \int_\E \big(\int_X M(y,z) u(t,z) \, dz  -\int_X M(z,y) \,dz \, u(t,y)  \big)  
\end{align*}
 The combination of the three mechanisms is illustrated in Figure \ref{figure 0} (right panel).

Taking the limit $\Delta t$ going to zero, the partial differential equation modelling the mechanisms of advection, growth and mutations is then written (for any dimension $n$) as
\begin{equation*}
\partial_t u(t,y)+ \nabla \cdot \left(f(y)u(t,y)\right)=\left(r(y)-d(y)\rho(t) \right)u(t,y) + \int_{\R^n}{M(y,z)u(t,z)dz}-\int_{\R^n}{M(z,y)dz}\,u(t,y),
\end{equation*}
which must be complemented with an initial condition describing the density of cells at time $0$.

\subsection{A Cell Population Balance Model recapitulates key dynamical aspects of EMT}

\begin{figure}[h]
    \centering
    \includegraphics[width=\textwidth]{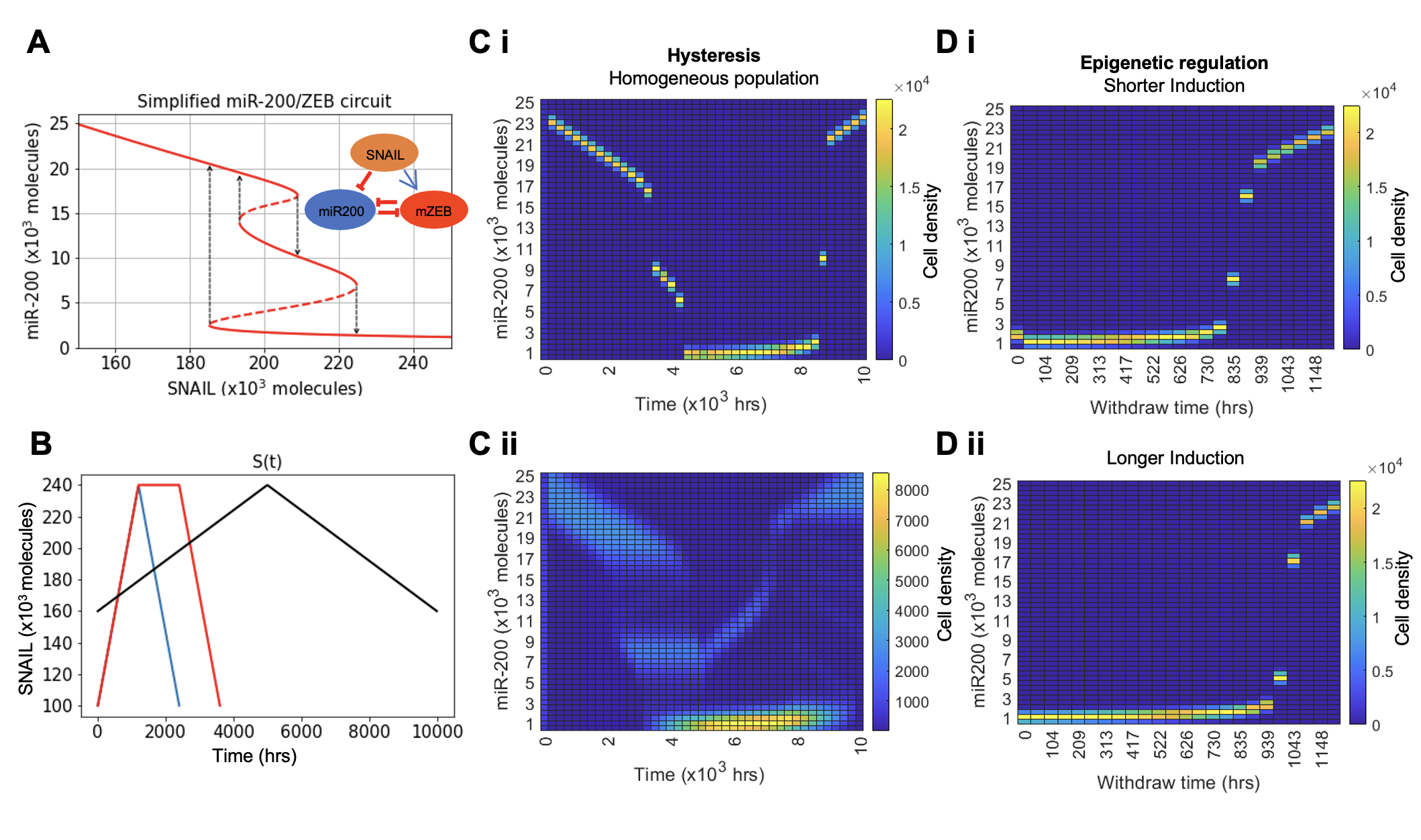}
    \caption{Proposed cell population balance model reproduces hysteresis and epigenetic regulation of state transitions at the population levels. A) EMT gene regulatory network (inset) and the bifurcation of cell states resulting from the network dynamics with increasing input signal (SNAIL) levels, B) Variation in the external input (SNAIL) levels with time to capture hysteresis and epigenetic regulation of cell states in C and D. C) Hysteresis (non-symmetric trajectories) in cell state transition during one cycle of EMT and MET by varying SNAIL levels (black curve) as shown in Figure 2B while considering – i) Homogeneous population and ii) Heterogeneous population. D) Dynamics of miR200 levels in homogeneous population during withdrawal of SNAIL levels post short and long-term EMT induction (blue and red curve in Figure 2B, respectively). Growth of cells and state transition because of biochemical noise are neglected in these simulation results.}
    \label{figure 2}
\end{figure}

\textit{In vitro} and \textit{in silico} studies on EMP have demonstrated two specific dynamic phenomena as cancer cells undergo one cycle of EMT and MET:
\begin{enumerate}
    \item Asymmetry in EMT and MET trajectories (hysteretic behaviour) of cell states. 
    \item Delayed MET with an increasing duration of EMT inducer treatment due to epigenetic (chromatin-based) stabilisation of M and hybrid E/M states.
\end{enumerate}
%– 1) Asymmetry in EMT and MET trajectories (hysteretic behaviour) of cell states.%, and 2) Delayed MET with an increasing duration of EMT inducer treatment due to epigenetic (chromatin-based) stabilization of M and hybrid E/M states. 
%For our further analyses, we choose a minimal EMT regulatory network with canonical epithelial (microRNA-200 (miR200)) and mesenchymal (ZEB) players that mutually inhibit each other via transcriptional and translational regulation (Figure \ref{figure 2}) . An EMT-inducing transcription factor SNAIL that activates ZEB and inhibits miR-200 represents the cumulative effect of the extracellular environment~\cite{lu2013microrna}. 
%As a starting point, we consider a simplification of the \textit{Core EMT Network} developed in \cite{lu2013microrna},
To properly define the advection function $f$ underlying state regulation, we chose a minimal EMT regulatory network with canonical epithelial (microRNA-200 (miR200)) and mesenchymal (mRNA ZEB) players that mutually inhibit each other via transcriptional and translational regulation. An EMT-inducing transcription factor SNAIL that activates ZEB and inhibits miR-200 represents the cumulative effect of upstream signalling pathways (Figure \ref{figure 2} A) \cite{lu2013microrna}. The bifurcation diagram depicts the different possible stable states, each characterised by a specific range of miR200 levels (solid lines) for increasing levels of SNAIL, resulting from the network dynamics. As a cell undergoes EMT (i.e. SNAIL levels increase), it switches from high to intermediate to low levels of miR200 which corresponds to E, E/M and M state respectively. However, during MET, the cell switches directly from low (M) to high (E) miR200 levels, thus displaying hysteresis.

We confirmed that the Cell Population Balance model developed here captures hysteresis, upon neglecting cell growth and transition (biochemical noise) and using SNAIL dynamics shown in black curve (Figure \ref{figure 2} B, see Methods Section for formalism). We simulated the dynamics for homogenous and heterogenous cell population with respect to their distribution of SNAIL levels. For a homogeneous cell population (Figure \ref{figure 2} C,i), we saw that the cells reside in three distinct miR200 states (high, intermediate, and low) for time $0$ to $5000$ hours (increasing SNAIL levels), but made a quick transition from low to high miR200 levels for time $5000$ to $10000$ hours without spending much time in intermediate state (decreasing SNAIL levels). The intermediate miR200 levels seen during MET are to be understood as a sample timepoint where miR200 levels are responding to changes in SNAIL levels before settling to its equilibrium (high) state. Similar observation of hysteresis was made while considering a heterogeneous cell population (Figure \ref{figure 2} C, ii). Particularly, the cell distribution along ZEB and miR200 axis during MET shows that the transient intermediate miR200 peak in homogeneous population has turned into little dispersed transient peaks which were clearly distinct from the intermediate peak arising during EMT (Figure \ref{figure S1}). This observation recapitulates the experimental data on different partial states seen during EMT vs. during MET \cite{karacosta2019Mapping}. Also, in the heterogenous population case, some cells fail to complete a full EMT, rather undergo a partial transition and then return to an epithelial state upon reduced SNAIL levels (Figure \ref{figure 2} C, ii). 

The range of SNAIL values for which E, E/M and M states are stable (Figure \ref{figure 2} A) can be altered by epigenetic (chromatin-based) changes that can happen during long-term EMT induction, leading to a delayed MET~\cite{jain2022Population, jia2019possible}. Therefore, to capture this phenomenon within our modelling framework, we incorporated the phenomenological formalism to account for epigenetic changes in EMT regulatory network (equation \eqref{core regulation equation epigenetic regulation}, \cite{jain2023epigenetic}). Again, we neglected cell growth and transitions (biochemical noise) for this analysis. We consider two SNAIL levels dynamics to show the influence of epigenetic changes during EMT: Short-term induction (Figure \ref{figure 2} B blue curve), and  Long-term induction (Figure \ref{figure 2} B red curve). We saw that homogeneous cell population have delayed recovery time for long-term EMT induction than short-term induction, therefore, recapitulating the previous observation based on population’s average cell analysis \cite{jain2023epigenetic}. 

%{\color{red} Figure S represents the distribution of cell population along ZEB and miR200 axis while the stability of E, E/M and M states changes during EMT and MET.} 

%When sticking to the exact PDE analogue of the ODE \eqref{core regulation equation 2D}, i.e., the advection equation, we recover the hysteretic phenomenon at the cell population level. In this case, we insist that cell growth and epimutations are not taken into account. The SNAIL dynamics used to obtain hysteresis are shown by the blue curve (Figure \ref{figure 2} B, Methods Section for formalism). 
Overall, the developed Cell Population Balance Model hence captures the dynamical features associated with EMT/MET.

\subsection{Biochemical reaction noise, coupled with regulatory cell-state dynamics, shapes heterogeneity pattern of the cell population}

Each individual cell's state (characterised by levels of a set of specific biomolecules) can dynamically evolve due to stochasticity in biochemical reactions or cell division, thus causing heterogeneity. We focus on the latter case, that is we consider that stochastic cell transitions occur at cell division. 

\begin{figure}[h!]
    \centering
    \includegraphics[width=\textwidth]{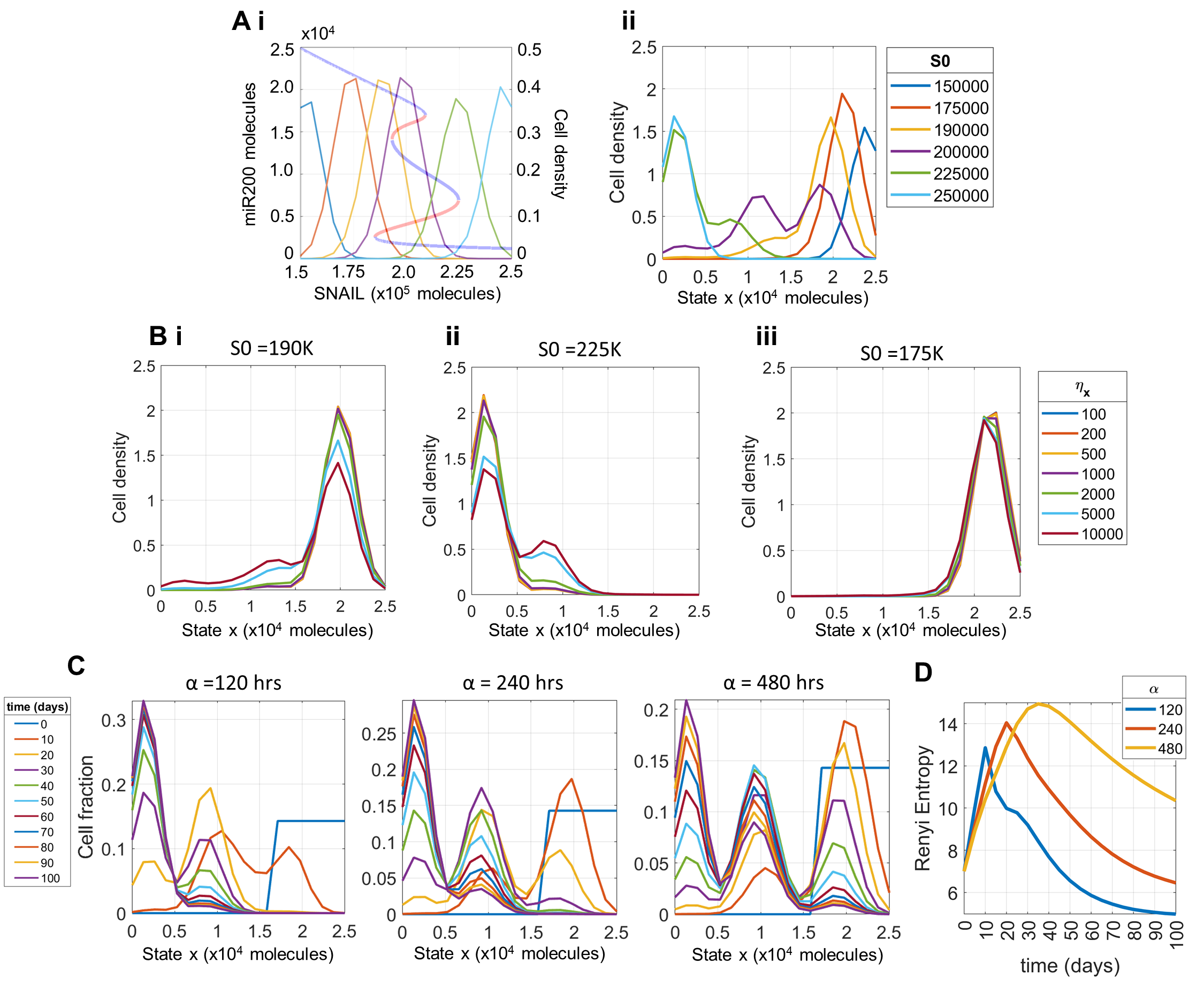}
    \caption{Simulations of the advection-selection-mutation model (PDE \eqref{advection-selection-mutation}) with reduced function $f_r$. Empirically reduced single cell state variable $x$, showing that the reduced system captures the full EMT network dynamical response; and the combined influence of biochemical noise with input signal dynamical response in modulating the phenotypic heterogeneity of the cell population. 
    Ai) Bifurcation diagram showing stable (blue) and unstable (red) levels of miR200 ($y$ axis on left) for increasing SNAIL levels; and distribution of SNAIL levels among the cells ($y$ axis on right) for different mean characteristics levels $S_0$.  Aii) Distribution of cells along the state $x$ of reduced EMT system for different distributions of SNAIL shown in Ai after a simulation time of 100 days starting from epithelial cells. B) Changes in population’s phenotypic distribution with increasing levels of epigenetic noise $\eta_x$ for SNAIL distributions with mean $S_0 = 190K$ (Bi) and 225K (Bii). For $S_0 = 175K$ (Bii), the population remains invariant to increasing noise.  C) Temporal changes in cell state distribution of the population for decreasing values of input signal SNAIL’s perturbation recovery rate (increasing values of the characteristic time ‘$\alpha$’). D) Time evolution of population’s heterogeneity (measured by Renyi entropy) for different $\alpha$ parameter values. Parameters used to generate plots, unless stated otherwise, are $\alpha= 120 hrs$, $\eta_x= 1000$, ini pop Epi (see Figure \ref{figure 4} A, i), time= 100 days, $S_0 = 225K$ molecules, and per-capita growth rate $r$ is constant across phenotypes, given by $r= 0.0182/ hr$.}
    \label{figure 3}
\end{figure}

To evaluate how these stochastic processes influence the population distribution of E, M and E/M states, we next observe cell state distribution as a result of stochastic cell-state transition in a regulatory network for a population where cell division rate is independent of cell-state (i.e. assuming EMT does not impact cell cycle)  (equation \eqref{advection-selection-mutation}) and is uniform for all three E, M and hybrid E/M subpopulations.
Because the levels of EMT-inducing signal SNAIL are also evolving due to biochemical reaction noise, its levels are distributed in the population around the mean environmental characteristics $S_0$ (Figure \ref{figure 3} A, i).
A cell defined by state $x, S$ gives birth to daughter cells defined by state $x',S'$. Our models captures how far $x'$ and $S'$ are from $x$, $S$ through standard deviations $\eta_x$ and $\eta_S$, respectively.

%\subsection{Focus on epimutations}    

To perform a comprehensive analysis of the impact of stochastic cell-state transition in a computationally efficient manner, we approximated the two-dimensional ODE (with variables miR200 and ZEB) in the EMT network by a one-dimensional ODE satisfied by a variable denoted $x$. The variable $x$ in the reduced system is built to be roughly equivalent to levels of miR200 in the full EMT network (more details about model reduction is presented in Appendix C). 

%We now turn to simulations for the full model~\eqref{biology advection-selection-mutation} where we investigate the effect of epimutations; growth functions $r$ and $d$ are kept constant.

Figure \ref{figure S2} A shows the shape of the function $x\mapsto f_r(x, S)$ underlying the reduced ODE for different values of $S$, and highlights the relative stability of possible cell-states for increasing levels of input S. The Methods section and Appendix \ref{annexes reduction} mention the empirical formulation of the reduced system and its optimal parameterisation to minimise the error in dynamical results obtained using complete vs reduced EMT system characteristics. We first established similarity between the dynamics of the system with two variables (full EMT) and that of the reduced system by comparing the distribution of miR200 levels for a given SNAIL distribution characteristic with the cell population distribution along the state ‘$x$’ (Figure \ref{figure 3} A and Figure \ref{figure S2} B). For example, with a SNAIL distribution of mean value ($S_0$) of 200K molecules, the bifurcation diagram depicts the possibility of cell population to be distributed in all three states (Figure \ref{figure 3} Ai). We observed that the asymptotic distribution of state variable $x$ from reduced system dynamics exhibits tri-modality, showing the co-existence of all three states, irrespective of initial population condition (Figure \ref{figure 3} A, ii – initial population: all cells as epithelial, Figure \ref{figure S2} B – initial population: all cells as hybrid or mesenchymal).  Similarly for SNAIL distribution with mean $S_0 \in \{190K, 225K, 150K, 250K\}$ molecules, we observed respective combinations of phenotypes as in bifurcation diagram of full EMT network – co-existence of E $\&$ E/M (bimodal), co-existence of hybrid E/M and M (bimodal), epithelial (unimodal) and mesenchymal (unimodal) (Fig \ref{figure 3} A), thus providing further evidence of faithfully representing EMT dynamics through the sole variable $x$. 
 
For multi-modal distributions of state variable $x$, the exact phenotypic composition depends on relative stability of the multiple co-existing phenotypes. For example, we see a reduced share of hybrid E/M (intermediate $x$ levels) cells for distribution of SNAIL levels that overlap significantly with those that have mean SNAIL levels corresponding to monostable E or monostable M regions ($S_0$ = 190K, 225K molecules respectively), especially for a reduced standard deviation $\eta_x$ of stochastic cell-state transition in state $x$ (Figure \ref{figure 3} Bi, ii). Similarly, for distribution of SNAIL levels with mean $S_0$ = 175K molecules, although both E and M states co-exist (Figure \ref{figure 3} A, i), the relative stability of the E state is much greater than that of the M state, thus disallowing  cells to make transition to the M state even at higher levels of $\eta_x$ (Figure \ref{figure 3} B, iii).

As mentioned previously, the distribution of SNAIL levels and correspondingly that of cell state $x$ in a population can be attributed to stochasticity in biochemical reactions. Another type of perturbation in cellular variables (here, $x$ and SNAIL) can arise when certain subpopulations are isolated and re-cultured independently. For instance, when the E, M and hybrid E/M prostate cancer subpopulations are segregated, they exhibit very different distributions after two weeks \cite{ruscetti2016HDAC}. Similarly, the segregated EpCAM-high and EpCAM-low subpopulations in breast cancer have varied recovery dynamics \cite{bhatia2019Interrogation}. Thus, in case of either internal (stochastic cell-state transition) or external (microenvironmental factors) perturbations, the rate at which the cellular variables recover towards the characteristic distribution can differ even though they eventually converge to the same equilibrium distribution.  Thus, we next modulated the rate of recovery of SNAIL levels to mimic the scenario of  extrinsic perturbation to the cell population by isolating distinct subpopulations and simulating (re-culturing) them independently.

The rate of recovery to perturbations in SNAIL levels is inversely proportional the parameter $\alpha$ in our model, which captures the characteristic time of convergence of SNAIL levels to its equilibrium $S_0$%. cells decorrelates to itself by $50\%$ 
~\cite{jain2023epigenetic, sigal2006Variability}.  %Starting from a mostly mesenchymal population (high values of $x$) the convergence to a mostly epithelial population slows down with increasing values of $\alpha$ (Figure \ref{figure 3} C).
For an extrinsic perturbation (e.g., enriched epithelial cells from a M cell majority population), we see that the time evolution of the population distribution slows down with increasing values of $\alpha$ (Figure \ref{figure 3} C). Furthermore, the slowed down dynamics increases cellular heterogeneity by causing the population to be distributed in all three states for a considerable interval of time, as quantified by Renyi entropy (Figure \ref{figure 3} D). %defined for a positive function $u\in L^1(\R^d)$ as $-\int_{\R^d}{\tilde u(x) \ln(\tilde u(x))dx}$, with $\tilde{u}:= u/\lVert u \rVert_{L^1}$.  
In the example shown, population heterogeneity first increases as the population shifts from a majority of epithelial cells to being more uniformly distributed among the three phenotypes in the intermediate time points, and then decreases as the population turns to a majority of M cells. Similar observations can be made for other combinations of initial condition and mean $S_0$ values of SNAIL distribution (Figure \ref{figure S3}, Figure \ref{figure S4}). 

Overall, the interplay between deterministic and stochastic dynamics of cellular %\linebreak 
biomolecules shapes the population heterogeneity. This is done both by distributing cells in all plausible states permitted by the underlying regulatory network dynamics, and by slowing down the kinetics of cells towards equilibrium when perturbed by the external signal $S$.  

\subsection{Difference in phenotypic growth rates reduces E-M heterogeneity, which could be recovered by increasing biochemical noise levels}

\begin{figure}[H]
    \centering
    \includegraphics[width=\textwidth]{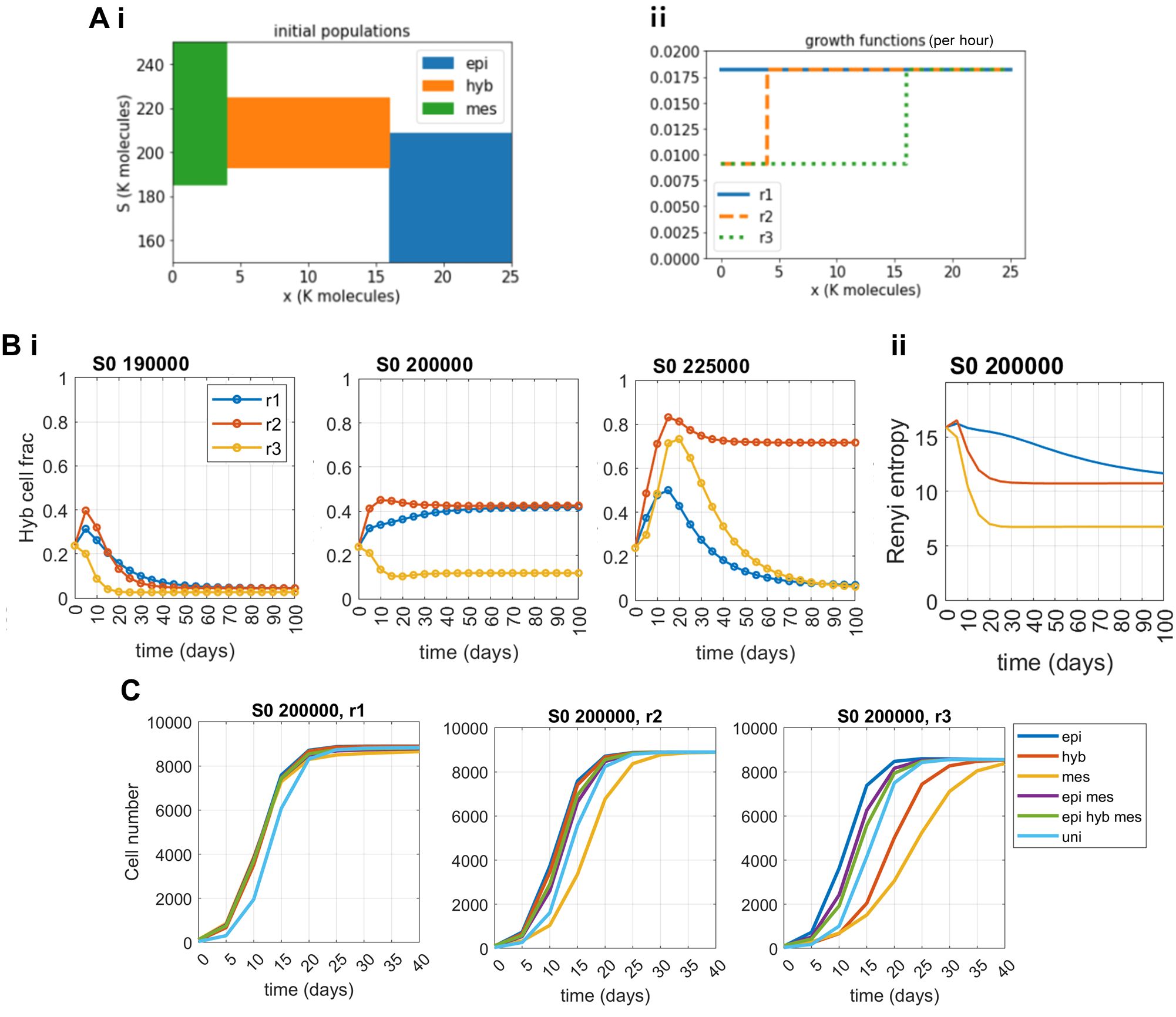}
    \caption{Effect of growth rate differences among E, hybrid E/M and M phenotypes on population’s heterogeneity and overall growth dynamics. Ai) Support of the different initial conditions: In each case, the initial condition is uniformly distributed on its support, and such that the total initial population equals 100 cells.  Aii) Profile of different growth functions, expressed per hour): ‘r1’: All three phenotypes divide at same rate; ‘r2’: E and E/M divide at equal rates, while M divide at half the rate of E cells; and ‘r3’: Both E/M and M divide at equal but half the rate of E cells. Bi) Temporal changes in hybrid cell fraction in the population for different growth scenarios among phenotypes.  Bii) Time evolution of population’s heterogeneity (measured by Renyi entropy). C) Population growth dynamics for different combinations of growth scenarios and initial condition; `epi hyb mes' corresponds to an initial condition uniformly distributed on the three colored domains of A1), and `uni' to a uniform population on the whole rectangle $[0, 25K]\times [150K, 250K]$. For panel~B, the initial condition is uniformly distributed in all states. The input SNAIL mean ($S_0$) levels used are mentioned for all the individual plots. Other parameters used to generate plots are $\alpha= 120 hrs$, $\eta_x= 1000$, and the per-capita growth rate (r) of the epithelial phenotype is $0.0182/ hr$.}
    \label{figure 4}
\end{figure}

We demonstrate how an interplay between deterministic and stochastic effects on the cell state $x$, and input SNAIL, can impact the E-M heterogeneity patterns. 
So far, we considered all three phenotypic states (E, M and hybrid E/M) to divide at equal rates. To observe the additional influence of growth rate differences on E-M heterogeneity we considered three possible scenarios – a) case ‘r1’: all three phenotypes divide at the same rate, b) case ‘r2’:  both E and hybrid E/M cells divide at equal rates, while M cells divide at half the rate of E cells; and c) case ‘r3’: both E/M and M divide at equal rates, which is half the rate of division of E cells. In practice, this means considering three different possible piecewise-constant functions $r$, as illustrated in Figure \ref{figure 4} A ii. Across all these cases, the E state divides at either an equal or a faster rate than hybrid E/M and/or M cells. This constraint recapitulates the current experimental understanding that EMT may suppress cell cycle to varying extents, thus reducing the division rate of hybrid E/M and/or M cells \cite{lovisa2015Epithelial, vega2004snail}. 

First, we investigate how the phenotypic composition of the population changes with different growth rate scenarios (Figure \ref{figure 4}, Figure \ref{figure S5}). For uniformly distributed cells in all three E (epi), hybrid E/M (hyb), and M (mes) states (Figure \ref{figure 4} A, i) and SNAIL distribution with mean $S_0 = 190K$ or 200K molecules that predominantly enables an E state with/without hybrid E/M state, the reduced growth of M cells has very slight effects on phenotypic composition over the time course, as expected (Figure \ref{figure 4} B, i). The initial peak in hybrid cell fraction for the ‘r2’ growth scenario is the combined effect of M to hybrid E/M state transition and a relatively higher growth rate of hybrid E/M cells. However, when hybrid E/M cells also have reduced proliferation (growth scenario ‘r3’), we see a lasting change in the phenotypic composition as E cells become dominant because of higher division frequency (Figure \ref{figure 4} B, i - $S_0 = 200K$). For the input SNAIL distribution with a mean value $S_0$ of 225K molecules that majorly supports hybrid E/M and M phenotypes, in the case ‘r2’, the growth advantage provided to hybrid E/M cells enables their dominance in the population on the long run, when compared to the growth scenarios of ‘r1’ or ‘r3’ where both E/M and M cells proliferate at equal rates (Figure \ref{figure 4} Bi, $S_0 = 225K$). The initial peaks in hybrid E/M fractions are combined effects of E to hybrid E/M transitions with either growth similarity or advantage of hybrid E/M cells over M cells. 

The overall change in phenotypic composition can be calculated using Renyi entropy as a heterogeneity score (Figure \ref{figure 4} B, ii). Although, growth scenarios ‘r1’ and ‘r2’ have the same phenotypic composition and an equal heterogeneity score eventually, the growth scenario ‘r1’ shows a much smoother change in heterogeneity values from the initial levels because of all three phenotypes being equally proliferative. Next, we look at the effects of increasing level of epigenetic noise level $\eta_x$ in cell state $x$ on phenotypic composition laid down by growth rate differences (Figure \ref{figure S6} A-D). For $S_0 = 190K$ and 200K, where hybrid E/M state is less dominant than the E state (Figure \ref{figure 4} B), increasing the noise levels (from $\eta_x = 1000$ to $\eta_x = 5000$) in state~$x$ causes more cell-state transitions, raising the frequency of hybrid E/M phenotype in population for all growth scenarios (Figure \ref{figure S6} B). However, as for $S_0 = 225K$ and growth scenario ‘r2’ where hybrid phenotype is the dominant state, increasing the noise $\eta_x$ level (from $\eta_x = 1000$ to $\eta_x = 5000$) raises M fractions in the population even though M cells are dividing slowly (Figure \ref{figure 4} B,i; Figure \ref{figure S6} B-C). Overall, we observe an increase in population heterogeneity with higher noise levels in state variable ‘$x$’, irrespective of the growth scenario (Figure \ref{figure 4} Bii, Figure \ref{figure S6} D).

After observing changes in phenotypic composition for different growth scenarios, we move on to see how total cell population grows for combinations of initial conditions and growth scenarios. We consider six different conditions – isolated E, hybrid E/M, and M population, uniform mixture of either E and M or E, E/M and M cells, and uniform distribution of cells in all possible cell states $x$ and input SNAIL levels. 

When all the phenotypic states are dividing at equal rates, the total number of cells does not vary across different initial conditions. However, with M dividing slower than E and hybrid E/M cells (growth scenario~‘r2’), we see that an initially mesenchymal population  has slower population growth compared to other initial conditions. Similarly, with both hybrid E/M and M cells dividing slower than E (growth scenario~‘r3’), initial conditions having component of either E/M or M cells divide slower than isolated (pure) E cells. Also, as transition from hybrid E/M to E is much more probable than M to E transitions, presence of hybrid cells in the populations increases the overall growth rate (Figure \ref{figure 4} C compare ‘hyb vs mes’, and ‘epi mes’ vs ‘epi hyb mes’ initial conditions for $S_0= 200K$, ‘r3’ growth scenario). Further, increasing the level of stochastic noise (from $\eta_x = 1000$ to $\eta_x = 5000$) causes more state transitions rendering lesser variability in the growth dynamics by quickly equalising the effect of differences in the initial conditions (Figure \ref{figure S6} E).

\subsection{Heterogeneity-dependent growth explains the faster tumour growth with highly heterogeneous parental and its subclones along E-M phenotypic axis}

Experimental data from orthotopic implantation of SUM149 cells and its subclones with varying degree of E, M and hybrid E/M heterogeneity levels suggest that the parental cell line and subclones with high levels of E-M heterogeneity enabled the fastest tumour growth in mice \cite{brown2022Phenotypic}. Further, despite starting from varying levels of E-M heterogeneity in the orthotopic injected cell populations, all tumours have relatively higher E-M heterogeneity levels when mice were sacrificed \cite{brown2022Phenotypic}. These observations together indicated a  plausible relationship between tumour growth rate and its heterogeneity. Thus, to explain these experimental observations, we assumed in our model formalism that the population growth rate depends on its heterogeneity levels (measured by differential entropy, a continuous equivalent of Renyi entropy). We specifically considered two functional relations between population growth rate and heterogeneity – 1) Linear relation, and 2) Sigmoidal relation (using Hill function with different threshold levels) (Figure 5A). 

\begin{figure}[H]
    \centering
    \includegraphics[width=\textwidth]{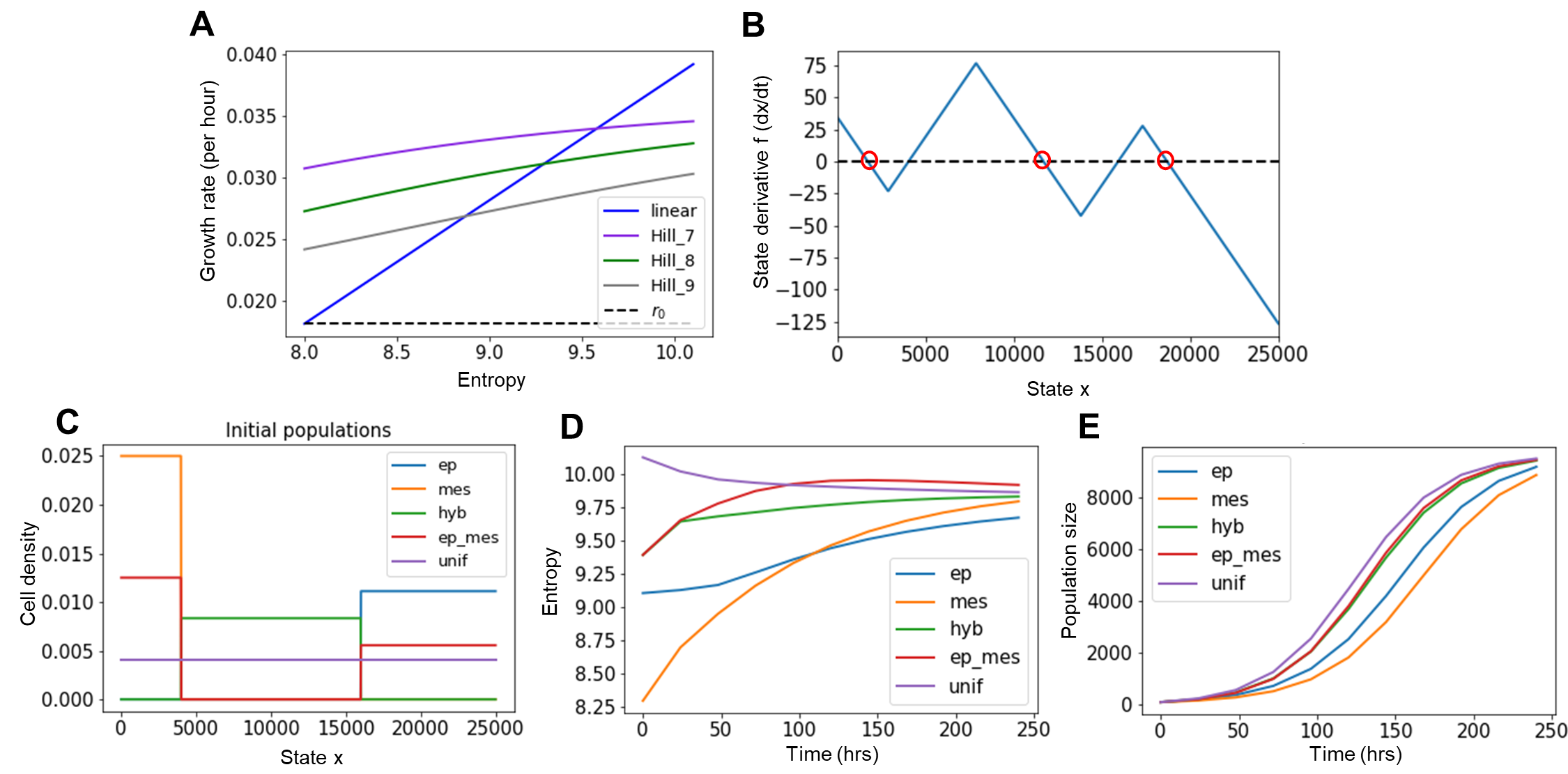}
    \caption{Heterogeneity-dependent growth of cells captures the in-vivo tumour growth dynamics observed in Figure 3A of \cite{brown2022Phenotypic}. A) Cellular division/growth rate as a function of population heterogeneity. Two functional dependences are considered – 1) Hills relation with increasing threshold value, and 2) Linear relation. More information about the functional relation is present in the Methods section.
    B) Reduced EMT state function~$f_r$ (underlying the reduced ODE) with three stable states (encircled) for constant S level of $200K$. C) Different starting population to simulate population dynamics.   D) Time evolution of the population heterogeneity (entropy) for different initial distributions of cells in E, E/M and M phenotypes when growth rates of cell have linear relation with entropy (linear in Panel A). E) Time evolution of total population size for simulation in panel D. Uniform (unif), Epithelial and Mesenchymal (ep\textunderscore mes), and Hybrid (hyb) starting populations grow faster because of their highly heterogeneous population at all time points. 
}
    \label{figure 5}
\end{figure}

To computationally track the population dynamics, we kept the S levels of cells to be invariant such that the stability landscape along cell state ‘x’ remained the same for all simulation times (Figure \ref{figure 5} B). We started our simulations with five different populations (distributed along state x) – Pure E (ep), pure M (mes) and pure hybrid E/M (hyb) populations, mixed E and M population (ep\_mes), and uniformly distributed population in E, M and E/M states (Figure \ref{figure 5} C). the difference in the heterogeneity among each starting population can be seen at zero-hour time point in Figure 5D. Further, as the cells divide and undergo cell-state transitions, the population heterogeneity value changes with time. All five different initial population distributions asymptotically reach saturating levels of heterogeneity, which is jointly determined by the stability landscape (Figure \ref{figure 5} B) and noise levels~($\eta_x$). Given that the population growth rate depends on its heterogeneity, we noticed that populations with higher E-M heterogeneity to start with (hyb, ep\_mes, and unif; Figure \ref{figure 5}D) show faster growth than other population with lower levels of E-M heterogeneity at the start (ep, and mes; Figure \ref{figure 5} D, E). The differences between the growth curves corresponding to different initial populations are much more pronounced for a linear than sigmoidal relationship between growth and heterogeneity because of its large variance in growth rates within the variability range of population’s heterogeneity (Figure \ref{figure 5} A, D, E and Figure \ref{figure S7} ).  Overall, by assuming the population growth rate to be dependent on its heterogeneity, we could recapitulate qualitatively the experimentally observed differences in the tumour growth dynamics in mice~\cite{brown2022Phenotypic}.

%\appendix

%\input{Annexes}

\section{Discussion}

Several intracellular and intercellular regulatory and stochastic process shape heterogeneity within a population. While several intracellular processes, such as non-linear transcriptional regulation, chromatin-based epigenetic regulation and microRNA-mRNA binding and complex degradation can result into the existence of distinct gene expression states, the stochastic intracellular processes such as transcriptional bursting, asymmetric cell division, and cell-to-cell communication lead the cells to switch from one gene expression pattern to the other. Experimental data capturing temporal changes in population level heterogeneity while profiling cells for transcriptome and epigenome is only recent in the context of EMT. However, significant efforts have been made in multi-scale mathematical modelling of a cell population that is growing, dividing and changing its phenotypic distribution with intracellular state dynamic based on one or more regulatory processes. Here, we contribute to this rich multi-scale population modelling literature by developing a framework allowing us to study how E-M population heterogeneity is regulated by – 1) the regulatory and stochastic intracellular processes and, 2) heterogeneity of growth rates among distinct subpopulations.

Our analysis is based on a minimalistic three node EMT regulatory network with a characteristic phenotypic landscape (Figure \ref{figure 2} A) ~\cite{lu2013microrna}. Our choice of a sufficiently simple such network was based on the following criteria – a) making the analysis computationally tractable, and b) integrating multiple processes together – cell division, death, cell-state transition and intracellular regulatory dynamics. However, many more complex EMT regulatory networks have been identified over the past decade ~\cite{hari2022Landscape}. Further, the minimalistic EMT model considered here does not have a regulatory term for the input signal, so we assumed the input to have negative feedback onto itself and that its levels in the population is distributed around the mean $S_0$. The ‘$\alpha$’ and $S_0$ in input ‘S’ dynamics (Equation \eqref{eq f_S}) can be considered as the inverse of the birth rate and ratio of birth to death rates of S molecules, respectively. Therefore, the parameter ‘$\alpha$’ sets maximum rate at which the population reaches its higher mean $S_0$ levels in the event of perturbation of the population distribution (Figure \ref{figure 3} C). With each molecule having its own birth and death kinetics, we see that the cellular memory can span across generations by inheriting the levels of a particular protein, as seen experimentally (Figure  \ref{figure 3}, ~\cite{corre2014stochastic, nordick2022nonmodular} ).    

To further reduce the long computation of the full EMT network when simultaneously considering cell division, death and mutation terms, we reduced the existing two-variable one input signal EMT state dynamical model to single-variable one input signal EMT model. The parameters of the reduced one-dimensional state derivative function (Figure \ref{figure S2} A) were set by minimising the error of its resulting dynamics with the evolution of cell state while two-dimensional state derivative (full EMT model) was considered. %While minimising the above error we ignored the division, death and epimutation events in the cell states. 
Our approximation and parameterisation of the state derivative function are in line with efforts over the last decade to use experimental flow cytometry and cell counts data to either parameterise the cell population balance model or optimally choose the functional form of state dynamics along with other parameters that fit the data well \cite{hasenauer2011identification, hasenauer2012analysis, loos2018hierarchical}. Although, these studies use Maximum Likelihood approaches to parameterise the system since their aim is to fit data rather than approximate an already existing but more complex model.
%which is more detailed parameter search algorithm than our grid search in the available parameter ranges in this study.

To model stochasticity in the cell’s state resulting from, for example, transcriptional bursting or asymmetric cell division, we have generalised intracellular noise by considering a general mutation kernel, which we then assumed to be normally distributed about the current cell state by focusing on the case where stochasticity is mediated by cell division.
 % over each infinitesimal time duration. 
We note that some literature has shown the individual contribution of each of the stochastic cellular process on cellular heterogeneity \cite{mantzaris2007single, mantzaris2006stochastic, shu2011Bistability}, and the model framework developed here can be easily used to account for other stochastic factors.

Our framework encompasses considering a diffusion term as has been done elsewhere. Indeed, a second-order diffusion term is the PDE counterpart of adding a Brownian motion to the considered ODE as in~\cite{jain2023epigenetic}. The corresponding diffusion term can be recovered by a suitable scaling of our integral mutation term as explained in~\cite{degond1989weighted}, which we have not done in the present work. % By considering Gaussian distribution over each infinitesimal time duration, we take into account all the higher order terms in the Taylor series expansion of the Gaussian distribution around its mean, which makes the noise generalization different from just diffusion (second order) approximation,  
Additionally, by approximating the noise to be normally distributed %over each infinitesimal time duration 
gave us a handle on its extent (standard deviation), and therefore, enabled us to comment on the changes in population heterogeneity with increasing levels of noise (Figure \ref{figure 3} B).
The main equation studied in this article \eqref{advection-selection-mutation} is in line with  Cell Population Balance models usually employed for heterogeneous cell populations ~\cite{waldherr2018estimation, spetsieris2012single, spetsieris2009novel, schittler2013new}. Nevertheless, we have chosen a rather common logistic shape for the growth term, which writes `$(r(y)-d(y)\rho(t))u(t,y)$', with  $\rho(t)$ the total number of cells at time~$t$, rather than a simple linear term as it is often done. 
From a biological point of view, this choice 
reflects the capacity for the cell population to self-regulate its growth due to density constraints. 
From the mathematical perspective, it guarantees that the size of the population does not blow up, \textit{i.e.} remains bounded, and allows one to study how the population evolves in larger times. 

Regarding the numerical implementation, we opted for a scheme from the family of particle methods: these schemes are indeed known to be well-suited in the context of PDEs with advection and `non-local' terms, that is terms that involve the density of the cell population at all points. In our case, the selection term and the mutation term are both non-local terms, since the former depends on the population size $\rho$, and the latter is a convolution with the so-called `mutation function' $M$. %Note that these schemes can also be used for equations with a diffusive term \cite{degond1989weighted}. 

Convergence of particle methods has been proved under conditions that are satisfied in our setting ~\cite{alvarez2023particle}. Compared with other methods used for the same type of problems, such as finite element methods~\cite{mantzaris2001numerical3} or finite difference methods~\cite{mantzaris2001numerical1}, particle methods have several advantages. Firstly, they are easily adaptable upon modifying the model, which allows for greater flexibility in model design (in our case, useful when passing from the homogeneous description \eqref{advection equation} to the full model with growth and epimutations \eqref{advection-selection-mutation}, and then to the entropy-dependent growth equation \eqref{entropy-dependent equation}). Moreover, they are based on a Lagrangian description of the system, and do not require an underlying mesh. More precisely, the initial data is discretised by a set of points, whose positions in the state space are then made to evolve in time via the advection term: this allows to `follow' the cell population as it converges towards regions of higher concentration. However, in the presence of mutation terms, particle methods are not asymptotic-preserving schemes~\cite{chertock2017practical}, which means that the particle approximation \eqref{particle approximation } does not correctly approach the solution of the PDE for very large times. To avoid this problem, we carry out the regulation process at each time step, as described in the Methods Section. 

By drawing a simplistic relationship from the experimental data between in vivo growth and tumour heterogeneity, we were able to explain trends in the tumour growth dynamics as seen (Figure \ref{figure 5}). However, we understand that various other factors such as feedback loops formed by interactions of tumor cells with extra-cellular matrix (ECM) and/or other stromal cells can alter heterogeneity patterns as well. For example, cells undergoing EMT can secrete LOXL2 that increases collagen crosslinking in ECM, and the ECM density as a well as stiffness can induce EMT.    

Overall, we employed cell population balance modelling to analyse the combined effect of EMT regulatory dynamics with cell division and death, and stochastic cell state transition. The integration of these complex processes together was made possible with an efficient PDE numerical integration scheme recently developed by some of us~\cite{alvarez2023particle}. 

\section{Methods}
\subsection*{Introduction to phenotype-structured PDE models} 
The state of a given cell is described by a time-dependent vector $y(t)\in \R^n$: for a given $i \in\{1, \ldots d\}$, $y_i(t)$ represents the concentration of some protein $i$ inside the cytoplasm of the cell at time $t$. In the context of EMT, the $y_i$'s can represent the level of several EMT markers such as miR-200 ($\mu_{200}$), ZEB ($Z$) or SNAIL~($S$).  The time-evolution of the cell state, for a single cell, is modelled by an Ordinary Differential Equation (ODE) of the form 
\begin{align}
\label{ode}
    \dot y(t)=f(t,y(t)),
\end{align}
where $\dot y$ denotes the derivative with respect to time, and $f : \R_+\times \R^n \rightarrow \R^n$ is a function describing the interactions between different molecules, which will be called \textit{advection function} throughout. Such an ODE is solved once complemented with an initial condition $y(0) = y^0$ where $y^0 \in \R^n$.

A population composed of many cells can be described by means of a density function $u(t,y)$, which represents the number of cells of state $y \in \R^n$ at time $t \geq 0$. In other words, the number of cells whose state lies in some set $\mathcal{E}\subset \R^n$ at a given time $t$ is given by $\int_{\mathcal{E}}{u(t,y)dy}$. The time-evolution of $u(t,y)$ associated to~\eqref{ode} is then given by a PDE, the so-called~\textit{advection equation}, namely 
\begin{align}
    \partial_t u(t,y)+\nabla \cdot \left( f(t,y) u(t,y) \right)=0, 
    \label{advection equation}
\end{align}
where $\partial_t$ denotes partial derivation with respect to time, and $\nabla \cdot$ denotes the (partial) divergence operator with respect to the variable $y$.
Such a PDE is solved once complemented with an initial condition $u(0, \cdot) = u^0$  with $u^0 : \R^n \to \R$.

\subsection*{Simulating hysteresis and epigenetic regulation}
We begin by considering the advection function associated with a minimal gene regulatory network developed in \cite{lu2013microrna}, which writes, 
\begin{equation}
\begin{cases}
\dot \mu_{200}&=g_{\mu_{200}}H_{Z,\mu_{200}}(Z)H_{S,\mu_{200}}(S)-g_{m_Z}H_{Z,m_Z}(Z)H_{S,m_Z}(S)Q(\mu_{200})-k_{\mu_{200}}\mu_{200}\\
\dot Z&=g_Zg_{m_Z}H_{Z,m_Z}(Z)H_{S,m_Z}(S)P(\mu_{200})-k_ZZ\\
\end{cases}.  
\label{core regulation equation 2D}
\end{equation}
We denote by $F = F(\mu_{200}, Z, S)$ the right-hand side of this ODE model, which accounts for interactions between a transcription factor ZEB (denoted $Z$), and a micro-RNA miR-200 (denoted $\mu_{200}$).  The variable $S$ represents a third molecule, SNAIL, which is seen in our case as an external signal characterising the extracellular environment. All the parameters of this model are given in Appendix \ref{Annexes Parameters ODE}.
Note that this model falls into the framework introduced in equation \eqref{ode}, considering that $y$ is the vector of miR-200 and ZEB concentrations $y = (\mu_{200},Z)$, and taking $f(t,y_1,y_2)= F(\mu_{200}, Z, S(t))$.

%Within this framework, the three cell states (E, E/M, and M) are characterised as follows:
%\begin{itemize}
%    \item The epithelial state corresponds to a high level of miR-200, and a low level of ZEB;
%    \item The mesenchymal state corresponds to a low level of  miR-200 and high level of ZEB;
%    \item The hybrid state corresponds to a medium level of miR-200 and ZEB.
%\end{itemize} 
%Depending on the value of $S$, system \eqref{core regulation equation 2D} exhibits monostability, bistability or tristability, as illustrated by the bifurcation diagram (Figure \ref{figure 1} A).

The bifurcation diagram displayed in Figure \ref{figure 2} A depicts the different possible stable states, each characterised by specific ranges of miR-200 levels (solid lines) for increasing levels of SNAIL, resulting from the network dynamics. As a cell undergoes EMT (i.e. as SNAIL levels increase), it switches from high to intermediate to low levels of miR-200, which corresponds to the E, E/M and M states respectively. However, during MET, the cell switches directly from low (M) to high (E) miR-200 levels without passing through the hybrid E/M state, thus displaying hysteresis. 
\medskip

\noindent
\textbf{Homogeneous population}, (Figure \ref{figure 2} C i). To reproduce the hysteretic behaviour of EMP, \textit{i.e.} asymmetry in EMT and MET trajectories, we first consider the advection equation~\eqref{advection equation},
where $y=(\mu_{200}, Z)$ is the two-dimensional vector of miR-200 and ZEB concentrations, the advection function is given by $f(t, \mu_{200}, Z)=F(\mu_{200}, Z, S(t))$, and $t \mapsto S(t)$ is the piecewise-linear function which connects the points $(0, 160K), (5000, 240K)$ and $(10000, 160K)$, as represented in Figure \ref{figure 2}B (black). 
\medskip

\noindent
\textbf{Heterogeneous population}, (Figure \ref{figure 2} C ii).
In order to account for heterogeneity within the population, and more specifically for the fact that the signal $S$ can be interpreted in a different way by each cell, we incorporate $S$ within the structure variable, which becomes $y=(\mu_{200}, Z, S)$. SNAIL level variation then impacts the advection term, which becomes $f(t,y)=(F(y), f_S(S))$, where $f_S$ is the step function corresponding to the derivative of the function $S$ introduced in the previous paragraph, \textit{i.e.} $f_S(S)=40$ if $S\in [0,5000)$, and $f_S(S)=-40$ if $S\in (5000, 10000]$.  As initial condition, we consider a population homogeneously distributed in the molecules ZEB and miR-200, but with heterogeneous levels of SNAIL distributed according to a Gaussian. In other words, we let $u^0(\mu_{200}, Z, S)=\frac{1}{\sigma}G(\frac{S-S_0}{\sigma})$, where $S_0=160K$, $\sigma=20K$ and $G$ denotes the Gaussian function. 
\medskip

\noindent
\textbf{Epigenetic regulation}, (Figure \ref{figure 2} D). Lastly, we  run simulations similar to those carried out for a homogeneous population, but with a modified advection function which allows to account for epigenetic regulation. We incorporate `$Z_0$' into the structure variable, which then writes $y=(\mu_{200}, Z, Z_0, S)$. This new parameter represents the ZEB threshold for inhibiting miR-200. The considered advection function is that associated to the ODE 
\begin{equation}
\begin{cases}
\dot \mu_{200}&=g_{\mu_{200}}H_{Z,\mu_{200}}( Z_0, Z)H_{S,\mu_{200}}(S)-g_{m_Z}H_{Z,m_Z}(Z)H_{S,m_Z}(S)Q(\mu_{200})-k_{\mu_{200}}\mu_{200}\\
\dot Z&=g_Zg_{m_Z}H_{Z,m_Z}(Z)H_{S,m_Z}(S)P(\mu_{200})-k_ZZ\\
\dot Z_0&=\frac{1}{\beta(t)} \left(Z^0_{\mu_{200}}-Z_0-\alpha Z  \right),
\end{cases}
\label{core regulation equation epigenetic regulation}
\end{equation}
whose parameters are detailed in Appendix \ref{Annexes Parameters ODE}. Denoting $F_e$ the right-hand side of this ODE, the corresponding PDE model~\eqref{advection equation} has advection function given by
%\begin{align*}
%    \partial_t u(t,y)+\nabla \cdot \left( f_e(t,y) u(t,y) \right)=0, 
%\end{align*}
$f_e(t, \mu_{200}, Z, Z_0, S)=F_e(\mu_{200}, Z, Z_0, S(t))$.
In a first simulation (Figure \ref{figure 2} D i),  $S$ is the piecewise-linear function which connects the points $(0, 100K), (1200, 240K)$ and $(2400, 100K)$, as represented in Figure~\ref{figure 2}~B (blue curve), and in a second simulation (Figure \ref{figure 2} D ii), the piecewise-linear function which connects the points $(0, 100K), (1200, 240K)$ and $(2400, 240K)$ and $(3600, 100K)$ (red curve in Figure \ref{figure 2} B).

The simulations for these three models are performed with a particle method detailed later on in this Methods section. 

\subsection*{Reducing the dimensions of the structure variable.} 
When incorporating growth and mutations into the model, computation times become prohibitive. To ease the burden, we use dimension reduction by further simplifying the advection function: we consider the state $y=(x, S)$ rather than $(\mu_{200}, Z, S)$, where $x$ is roughly equivalent to $\mu_{200}$, as explained below.

This requires to choose a function $f_r$ such that the dynamics of $x$ for any $S$ is given by $\dot x=f_r(x, S)$. Our main requirement in choosing this new function is to preserve the bifurcation diagram given in Figure~\ref{figure 2} A, which means that, for a given value of $S\in [150K, 250K]$, $f_r(\cdot, S)$ has the same number of zeros as $F(\cdot, S)$, and are such that $f_r(x, S)=0$ if and only if there exists $Z>0$ such that $F(x, Z, S)=0$. 

Since infinitely many functions satisfy this property, one must make further assumptions in order to select a suitable one. For simplicity and to avoid overfitting, we assume that, for any $S$, $f_r(\cdot, S)$ is piecewise linear (with one root per interval where it is linear) and that the rate of change is constant on each interval. Under these constraints, the function is defined up to a multiplicative constant which is chosen by minimising a suitably defined criterion, see Appendix \ref{annexes reduction}. 

\section*{Considering growth and epimutations.}
The full PDE model incorporating growth and epimutations writes
\begin{equation}
\begin{cases}
\partial_t u(t,y)+ \nabla \cdot \left(f(y)u(t,y)\right)=\left(r(y)-d(y)\rho(t) \right)u(t,y)\\
\hspace{1.4 cm} +\int_{\R^2}{M(y,z)u(t,z)dz}-\int_{\R^2}{M(z,y)dz}\,u(t,y)\\
\rho(t)=\int_{\R^2}{u(t,y)dy}
\end{cases},
\label{advection-selection-mutation}
\end{equation}
where:
\begin{itemize}
\item The structure variable $y\in \R^2$ is $y=(x, S)$. 
\item The advection function $f$ is given by $f(y)=f(x,S)=(f_r(x,S), f_S(S))$, where 
\begin{equation}
    f_S(S)=\delta \big(1-\frac{S}{S_0}\big),
\label{eq f_S}
\end{equation}
with $S_0\in [150K, 250K]$ corresponding to the mean of the SNAIL distribution, $\delta:=\frac{S_0 \ln(2)}{\alpha}$, and $\alpha>0$ representing the characteristic time of convergence of SNAIL to the mean $S_0$. 
\item The term $\left(r(y)-d(y)\rho(t) \right)u(t,y)$ represents the net growth of cells of state $y$, with two main contributions given by the intrinsic growth rate $r(y)$, and the death rate $d(y) \rho(t)$ proportional to the total population size $\rho(t)$. This corresponds to the so-called logistic model, accounting for the additional death rate due to competition for resources and space between cells. 
%\item $r(y)$ and $d(y)$ respectively represent the birth and death rate of a cell of state~$y$. 
In all our simulations, the death rate is considered to be independent of the cell state (\textit{i.e.} $d(y)=d\equiv 1.82 \times 10^{-7}$cell/hr, as in \cite{tripathi2020mechanism}).
\item The last two terms represent cell mutations, and can again be decomposed into two terms.
The term $\int_{\R^2}{M(y,z)u(t,z)dz}$ represents the mutation of cells of any type $z$ into cells of type $y$, occurring with a rate $M(y,z)$. The term $\int_{\R^2}{M(z,y)dz}\,u(t,y)$ represents the mutation of cells originally of type $y$ into cells of any other type $z$, with mutation rate $M(z,y)$.
The mutation function $M$ is taken to be $M(y, z)=r(z) P(y-z)$, meaning that mutations are considered to occur at cell division.
%$M(y, z)=\frac{r(z)}{\tau}P(y-z)$. 
Here, $P(y)=P(x,S):=\frac{1}{\eta_x \eta_S}G(\frac{x}{\eta_x})G(\frac{S}{\eta_S})$,
where $G$ is the Gaussian function. Variables $\eta_x$ and $\eta_S$ are the standard deviations for $x$ and $S$ respectively.%, and $\tau$ is a transition time characterising the frequency at which cells mutate, which is assumed to be constant equal to $1$ throughout. 

%($G(s)=\frac{1}{\sqrt{2\pi}}e^{-\frac{s^2}{2}}$).
\end{itemize}

\subsubsection*{Numerical method}

For numerical simulations, we start by normalising the model to work with the domain $[0,1]\times [0,1]$ rather than $[0, 25K] \times [150K, 250K]$. By denoting $A=25K$, $B=0$, $C=100K$, $D=150K$, we check that for all $t\geq 0$, $u(t, x,S)=\frac{1}{AC}u_{re}(\frac{x-B}{A}, \frac{S-D}{C})$, where $u_{re}$ is the solution of 
\begin{equation}
\begin{cases}
\partial_t u_{re}(t,y)+ \nabla \cdot \left(f_{re}(y)u_{re}(t,y)\right)=\left(r_{re}(y)-d_{re}(y)\rho(t) \right)u_{re}(t,y)\\
\hspace{1.4 cm} +\int_{\R^2}{M_{re}(y,z)u(t,z)dz}-\int_{\R^2}{M_{re}(z,y)dz}\,u(t,y)\\
\rho(t)=\int_{\R^2}{u_{re}(t,y)dy}
\end{cases},
\label{advection-selection-mutation_re}
\end{equation}
where for all $x, S \in \R$,  $f_{re}(x,S):=(\frac{1}{A}f_r(Ax+B,CS+D), \frac{1}{C}f_S(CS+D))$, $r_{re}(x,S):=r(Ax+B, CS+D )$, $d_{re}(x,S):=d(Ax+B, CS+D )$, $M_{re}(x,S, x',S'):=ACM(Ax+B,CS+D, Ax'+B, CS'+D)$, and $u_{re}^0(x,S):=ACu^0(Ax+B, CS+D)$. 

Thus, an approximation for $u_{re}$ provides an approximation for $u$. We apply a particle method in order to approximate $u_{re}$ at different time steps ($0<T_1<...<T_K$, specified in each figure), applying a particle method introduced in~\cite{alvarez2023particle} to deal with a category of models to which~\eqref{advection-selection-mutation} belongs. For a proof that the numerical scheme does successfully approximate the solutions of~\eqref{advection-selection-mutation}, we refer to~\cite{alvarez2023particle}, while an introduction to particle methods can be found in \cite{chertock2017practical}. 
%and more details on this method can be found in this paper, as well as 

We choose an integer parameter $N$ ($N=20$ in our simulations), and we denote $y_1^0, \ldots, y_{N^2}^0$, the points of the grid of size $N \times N$ on $[0,1]\times [0,1]$. 
For $i\in \{1,..., N^2\}$, we solve the ODE  
\begin{equation}
    \begin{cases}
        \dot y_i=f(y_i)\\
        \dot w_i=\nabla \cdot f (y_i)\, w_i\\
        \dot v_i= \left(r(y_i)-d(y_i)\sum\limits_{j=1}^{N^2}{v_j}\right)v_i+w_i\sum\limits_{j=1}^{N^2}{M(y_i, y_j)v_j}-\sum\limits_{j=1}^{N^2}{M(y_j, y_i)}v_i
    \end{cases}, 
    \label{particle approximation }
\end{equation}
with initial conditions $y_i(0)=y_i^0$, $w_i(0)=\frac{1}{N^2}$ and $v_i(0)=\frac{u_{re}^0(y_i^0)}{N^2}$, on the first time interval ($[0,T_1]$). The solution of this ODE is called \textit{Particle approximation of \eqref{advection-selection-mutation}}.
To solve \eqref{particle approximation } on $[0,T_1]$, we use the Python function solve\textunderscore ivp in module scipy.integrate, with the default solver which corresponds to an explicit Runge-Kutta method of order 5 \cite{dormand1980family}. 

We then use a \textit{regularisation process}, \textit{i.e.} we compute the sum
\begin{equation}
u^N_\varepsilon(T_1,x)=\sum\limits_{i=1}^{N^2}v_i(T_1)  \, G_{\varepsilon}\big(y-y_i(T_1)\big), 
\end{equation}
with $G_\varepsilon(x,S):=\frac{1}{\varepsilon^2} G(\frac{x}{\varepsilon})G(\frac{S}{\varepsilon})$, with $G$ the Gaussian function, and $\varepsilon=\left(\frac{1}{N^2}\right)^{\gamma}$, with $\gamma\in (0.5, 1)$. In all our simulations, we choose $\gamma=0.8$: this value has been chosen empirically by carrying out simulations in simple cases for which the behaviour of solutions is well known (for example in the absence of mutations), and by comparing with other values of $\gamma$. 
The points `$y$' at which we compute this sum are $y_1^0, \ldots, y_{N^2}^0$. 

We repeat the process for each time interval, taking as initial data the approximation calculated at the previous time step (\textit{i.e.}, $u^N(T_{k-1}, \cdot)$), to compute $u^N(T_{k}, \cdot)$.

\subsection*{Entropy-dependent growth. }

In Figure \ref{figure 5}, we consider a model for which the growth function $r$ depends on population heterogeneity. For simplicity, we assume that the SNAIL level is constant within the population, and does not change over time. Thus, we use the simplified advection function $f_r$ with $S\equiv 200 K $, which allows for the existence of the three phenotypes as shown on the bifurcation diagram (Figure \ref{figure 2} A). We measure heterogeneity within the population by means of the entropy 
$-\int_{\R}{\frac{u(t,y)}{\rho(t)} \ln ( \frac{u(t,y)}{\rho(t)})dy}$. The model writes 
\begin{align}
\begin{cases}
\partial_t u(t,x)+ \nabla \cdot \left(f(x)u(t,x)\right)=\left(r(x, E(t))-d(x, E(t))\rho(t) \right)u(t,x)\\
 +\int_{\R}{M(x,y, E(t))u(t,y)dy}-\int_{\R}{M(y,x, E(t))dy}\,u(t,x)\\
\rho(t)=\int_{\R}{u(t,y)dy}\\
E(t)=-\int_{\R}{\frac{u(t,y)}{\rho(t)} \ln \left( \frac{u(t,y)}{\rho(t)} \right)dy}
%\partial_t u(t,x)+\nabla \cdot (f(x)u(t,x))=\bigl(r(x)+E(t)-d(x)\rho(t)\bigr)u(t,x)+\frac{1}{\tau}\int_{\R^n}{G(y-x)u(t,y)dy}-\frac{1}{\tau}\int_{\R^n}{G(y-x)u(t,x)dy} \\
%\rho(t)=\int_{\R^n}{u(t,y)dy} \\
%E(t)=-\int_{\R^n}{\frac{u(t,y)}{\rho(t)} \ln %\left( \frac{u(t,y)}{\rho(t)} \right)dy}\\
%u(0,x)=u^0(x)
\end{cases}, 
\label{entropy-dependent equation}
\end{align}
with $f(x):=f_r(x, 200K)$, and $M(x,y)= \frac{1}{\eta_x}G\big(\frac{x-y}{\eta_x}\big)r(y,E)$ with $G$ the Gaussian function and $\eta_x=4000$, $d(x, E)=r(x,E)/10K$. We take four different values for the growth function `$r$': shifted Hill functions of the shape $r_0 \frac{\theta^6+2E^6}{\theta^6+E^6}$, with $\theta \in \{7, 8, 9\}$, $r_0=0.0182$ ($\theta = 9$ in Figure \ref{figure 5}), and a linear function $r_0+0.01(I-8)$. These four functions are represented in Figure \ref{figure 5} C. 

%def r(x,I):
%    return(1*(0.0182+0.01*(I-8)))

%$d(x,E)=r(x,E)/10K$,    $M(x,y)= 1./\eta_xG((x-y)/\eta_x)r(x,E))$, $\eta_x=4000$, where $G$ denote the Gaussian function. 

We use the same numerical scheme as for the previous equation (in one dimension), but with the following ODE
\begin{equation}
    \begin{cases}
    \dot x_i=f(x_i)\\
        \dot w_i=\nabla \cdot f (x_i)\, w_i\\  
        \dot v_i= \left(r(x_i, \bar E)-d(y_i, \bar E) \bar \rho\right)v_i+w_i\sum\limits_{j=1}^N{M(y_i, y_j, \bar E)v_j}-\sum\limits_{j=1}^N{M(y_j, y_i, \bar E)}v_i\\
        \bar \rho=\sum\limits_{j=1}^N{v_j}, \quad \bar E=-\sum\limits_{j=1}^N{\frac{v_j}{\bar \rho}\ln ( \frac{v_j}{w_j \bar \rho} )}
    \end{cases},
\end{equation}
with parameters $N=50$ and $\varepsilon=\left(\frac{1}{N}\right)^{0.8}$.

\paragraph{Funding.}
This work was supported by the Raman-Charpak Fellowship 2022 awarded to J.G., funding a two-month research stay at the Indian Institute of Science, Bangalore. M.K.J. was supported by Ramanujan Fellowship (SB/S2/RJN-049/2018) awarded by Science and Engineering Research Board (SERB), Department of Science and Technology (DST), Government of India. N.P.D. was supported by Emergence fellowship (S21JR31024) awarded by Sorbonne University.

\paragraph{Author Contributions.}
J.G. and P.J. performed research. M.K.J., C.P. and N.P.D. designed and supervised research. All authors contributed to data analysis and in writing and reviewing the manuscript.
\newpage

\appendix

\beginsupplement

\section{Supplementary figures}

\begin{figure}[H]
    \centering
    \includegraphics[width=\textwidth]{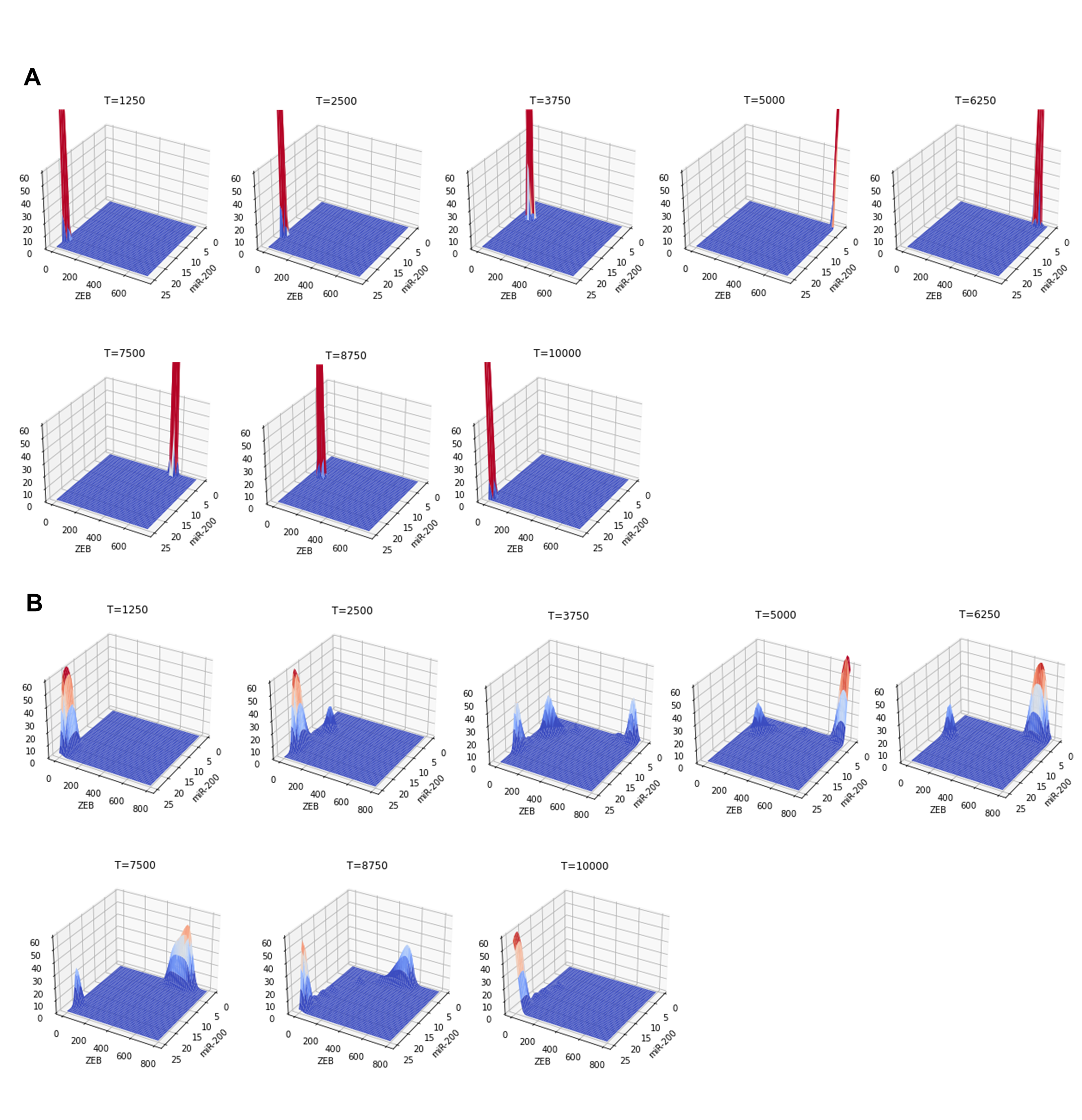}
    \caption{Hysteresis (non-symmetric transition) in cell density (z-axis) along the two cell state variables miR200 and ZEB during one cycle of EMT and MET caused by increasing and decreasing levels of input SNAIL levels (Figure \ref{figure 2} B blue curve) for A) homogeneous and B) heterogeneous cell population.}
    \label{figure S1}
\end{figure}

\begin{figure}[H]
    \centering
    \includegraphics[width=\textwidth]{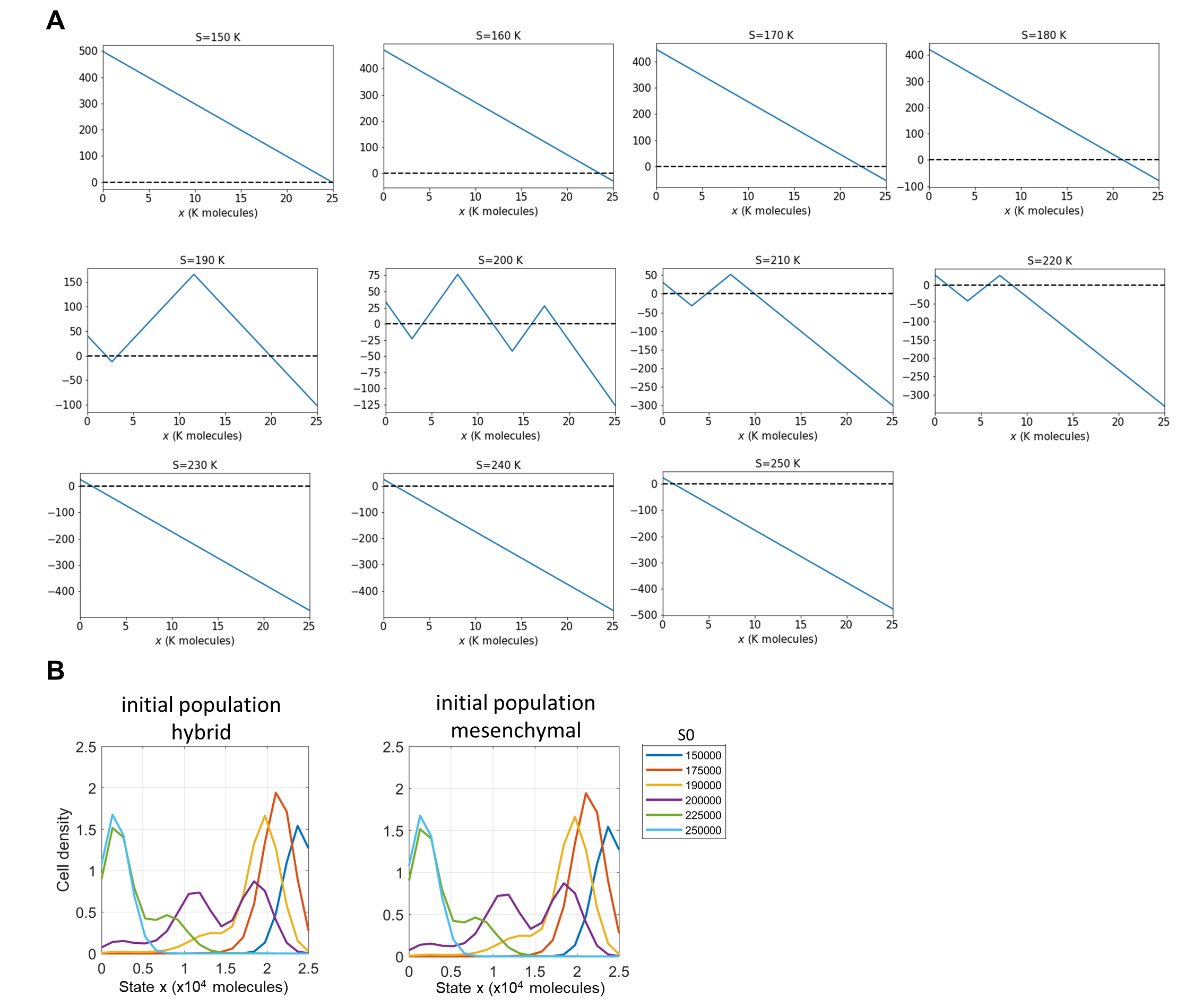}
    \caption{Reduced function $f_r$ underlying the ODE with single variable cell state $x$ and single input SNAIL; and distribution of the cell population for different distributions of SNAIL shown in Figure \ref{figure 2} Ai. A) Function $x \mapsto f_r(x,S)$, showing the existence of mono-, bi-, and tri-stable states for varying levels of SNAIL inputs $S$. B) Cell population distribution at the end point (100 days) of simulations started with hybrid and mesenchymal populations for increasing levels of input signal SNAIL’s mean ($S_0$) levels. Parameters used to generate the above plots, unless stated otherwise, are $\alpha = 120$ hrs, %$\tau =1$, 
    $\eta_x =1000$, ini pop Epi, time point 100 days, $S_0 = 200K$ molecules, and the per-capita growth rate (r) is constant $= 0.0182\backslash hr$.}
    \label{figure S2}
\end{figure}

\begin{figure}[H]
    \centering
    \includegraphics[width=\textwidth]{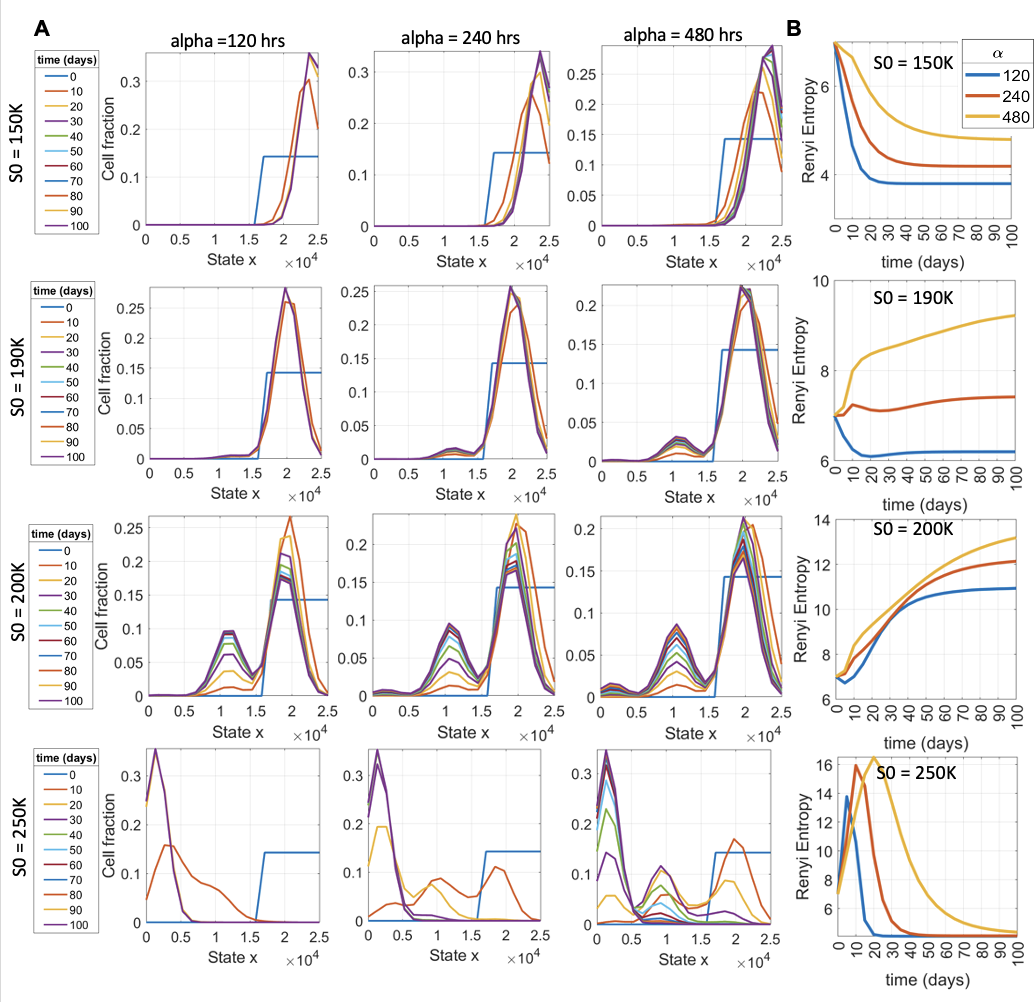}
    \caption{Population dynamics and changes in population heterogeneity (measured using Renyi Entropy) with time for several combinations of $S_0$ and $\alpha$ parameters. The intermediate $S_0$ values of 200K molecules where a cell can either attain a stable epithelial, hybrid, and mesenchymal, has highest heterogeneity score. Further, as we go towards smaller or larger values of $S_0$ where only epithelial, and mesenchymal states are possible, respectively, the population has least heterogeneity. And, by increasing the residence time of cells in a state $x$, cells get enough time to populate the cell state they reside in, and thereby, contribute to a significant fraction in the overall population. For the above plots, the initial condition population is epithelial, %$\tau = 1$ 
    $\eta_x = 1000$, and per-capita growth rate $r$ is constant $= 0.0182 / hr.$} %A time step is equal to 120 hrs. }
    \label{figure S3}
\end{figure}

\begin{figure}[H]
    \centering
    \includegraphics[width=\textwidth]{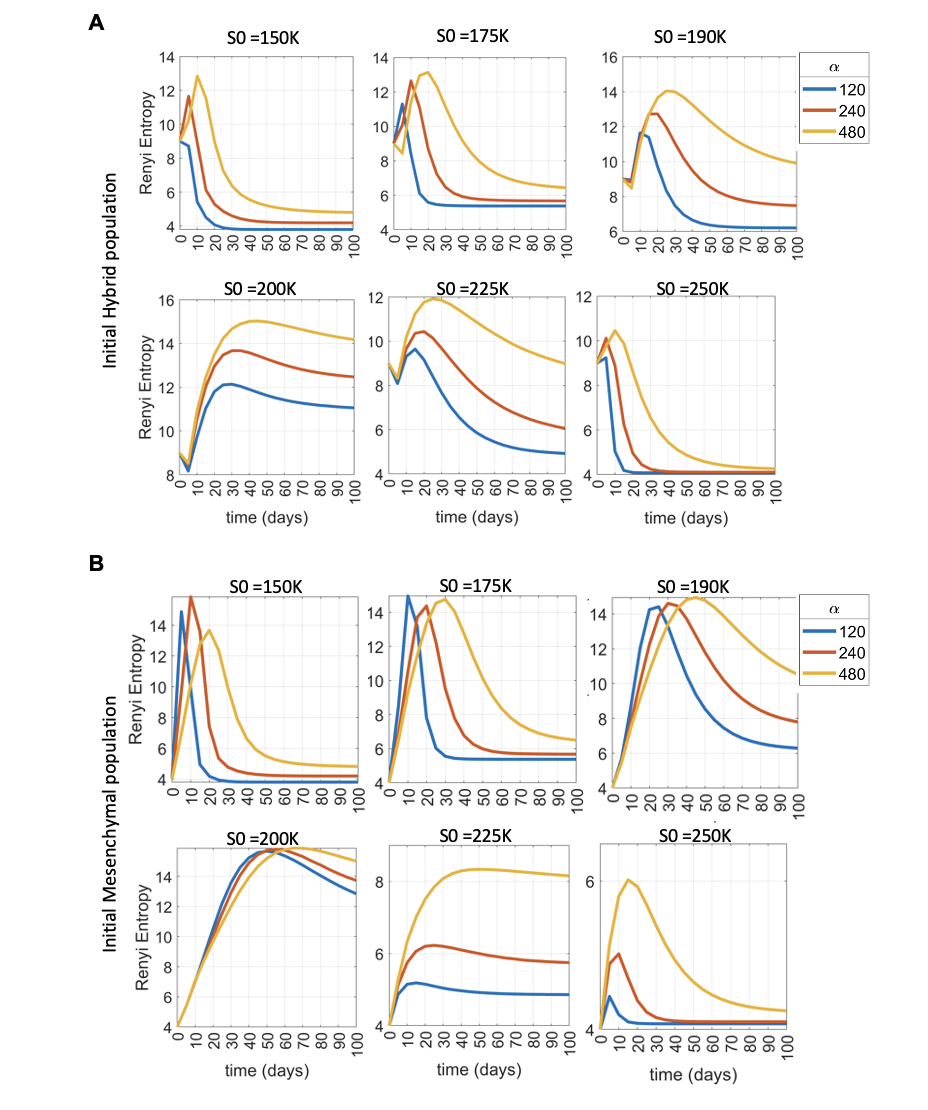}
    \caption{Changes in population heterogeneity (measured using Renyi Entropy) with time for several combinations of $S_0$ and $\alpha$ values while starting with a population of A) Hybrid cells, and B) Mesenchymal cells. For the above plots, %$\tau = 1$ 
    and $\eta_x = 1000$, and per-capita growth rate $r$ of all subpopulation = 0.0182/hr. %A time step is equal to 120 hrs.
    }
    \label{figure S4}
\end{figure}

\begin{figure}[H]
    \centering
    \includegraphics[width=\textwidth]{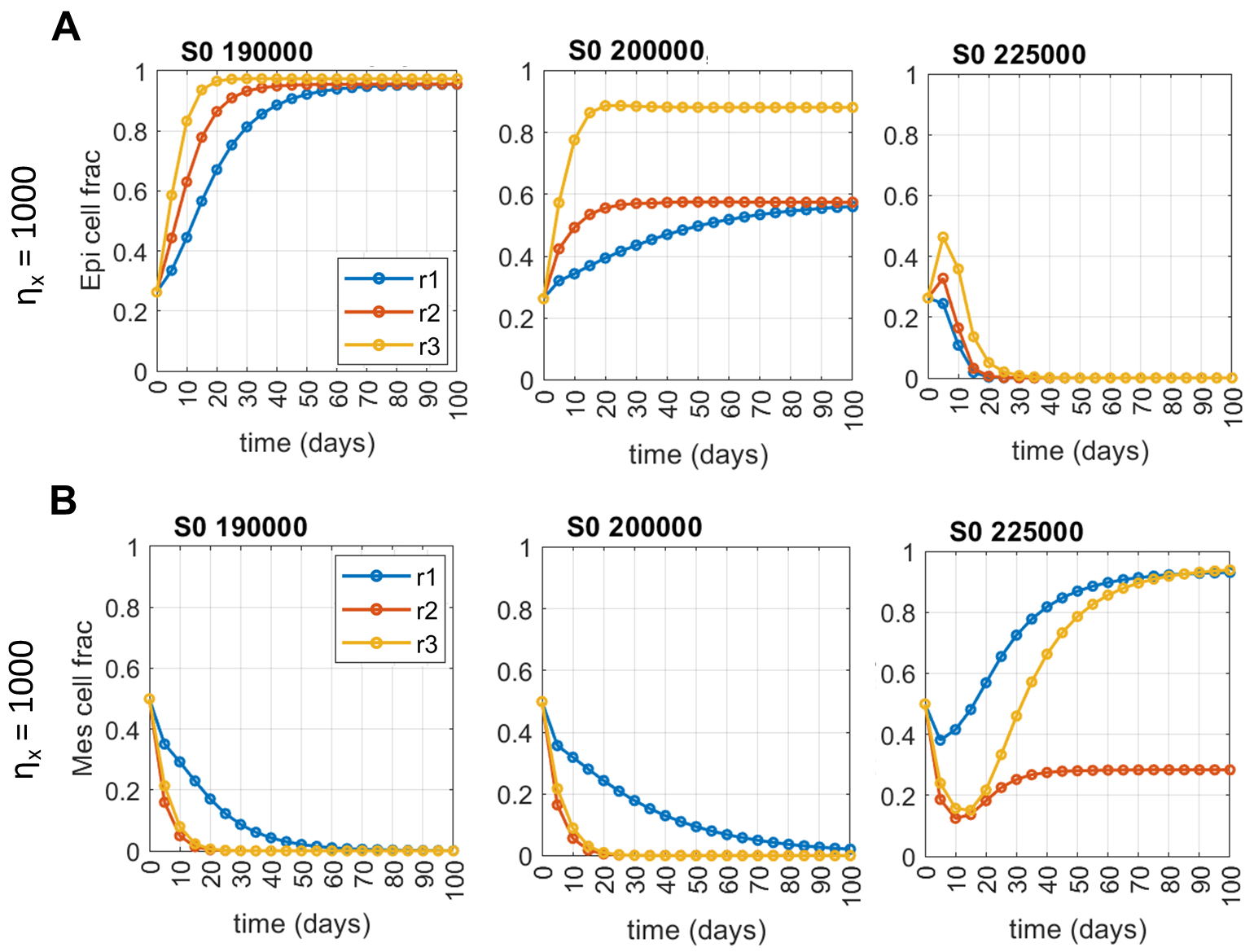}
    \caption{Effect of growth rate differences among E, hybrid E/M and M phenotypes on population’s heterogeneity. Temporal changes in E cell fraction (panel A) and M cell fraction (panel B) in the population for different growth scenarios among phenotypes – ‘r1’: All three phenotypes divide at same rate; ‘r2’: E and E/M divide at equal rates, while M divide at half the rate of E cells; and ‘r3’: Both E/M and M divide at equal but half the rate of E cells. For the above results, the initial condition is uniformly distributed in E, hybrid E/M and M state. The input SNAIL mean $S_0$ levels used are mentioned for all the individual plots. Other parameters used to generate plots are $\alpha= 120$ hrs, $\eta_x =1000$, and per-capita growth rate $r$ of epithelial phenotype is 0.0182/hr.}
    \label{figure S5}
\end{figure}

\begin{figure}[H]
    \centering
    \includegraphics[width=\textwidth]{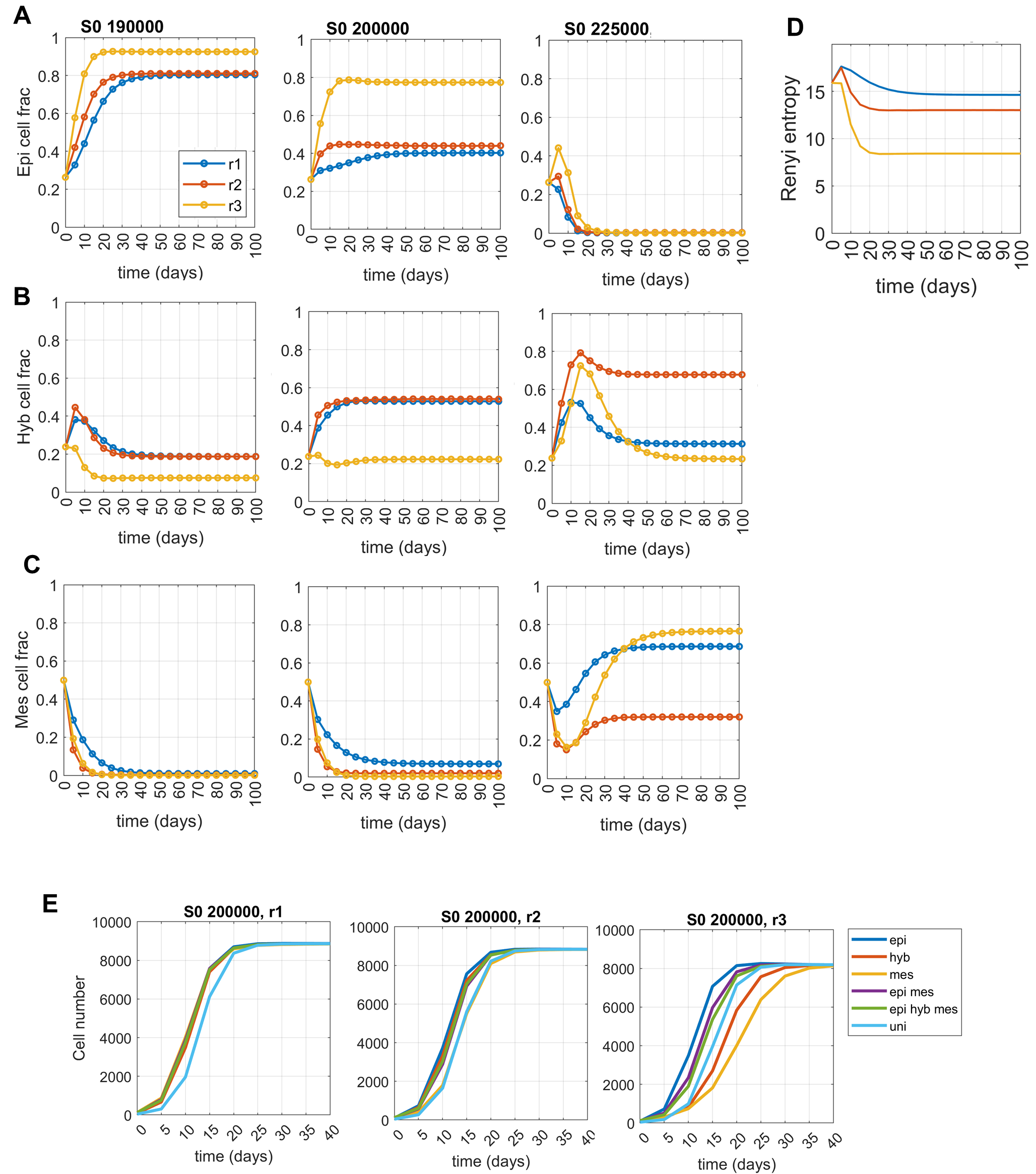}
    \caption{ Effect of growth rate differences among E, hybrid E/M and M phenotypes on population’s heterogeneity and overall growth dynamics. A-C) Temporal changes in E, hybrid E/M, and M cell fraction in the population for different growth scenarios among phenotypes – ‘r1’: All three phenotypes divide at same rate; ‘r2’: E and E/M divide at equal rates, while M divide at half the rate of E cells; and ‘r3’: Both E/M and M divide at equal but half the rate of E cells. D) Changes in population’s heterogeneity (measured by Renyi entropy) with time. E) Population growth dynamics for different combinations of growth scenarios and initial  conditions. For panel A-D, the initial condition is uniformly distributed in E, hybrid E/M and M state. The input SNAIL mean $S_0$ level used are mentioned for all the individual plots. Other parameters used to generate plots are  $\alpha= 120$ hrs, $\eta_x =5000$, and per-capita growth rate $r$ of epithelial phenotype is 0.0182/hr.
}
    \label{figure S6}
\end{figure}

\begin{figure}[H]
    \centering
    \includegraphics[width=\textwidth]{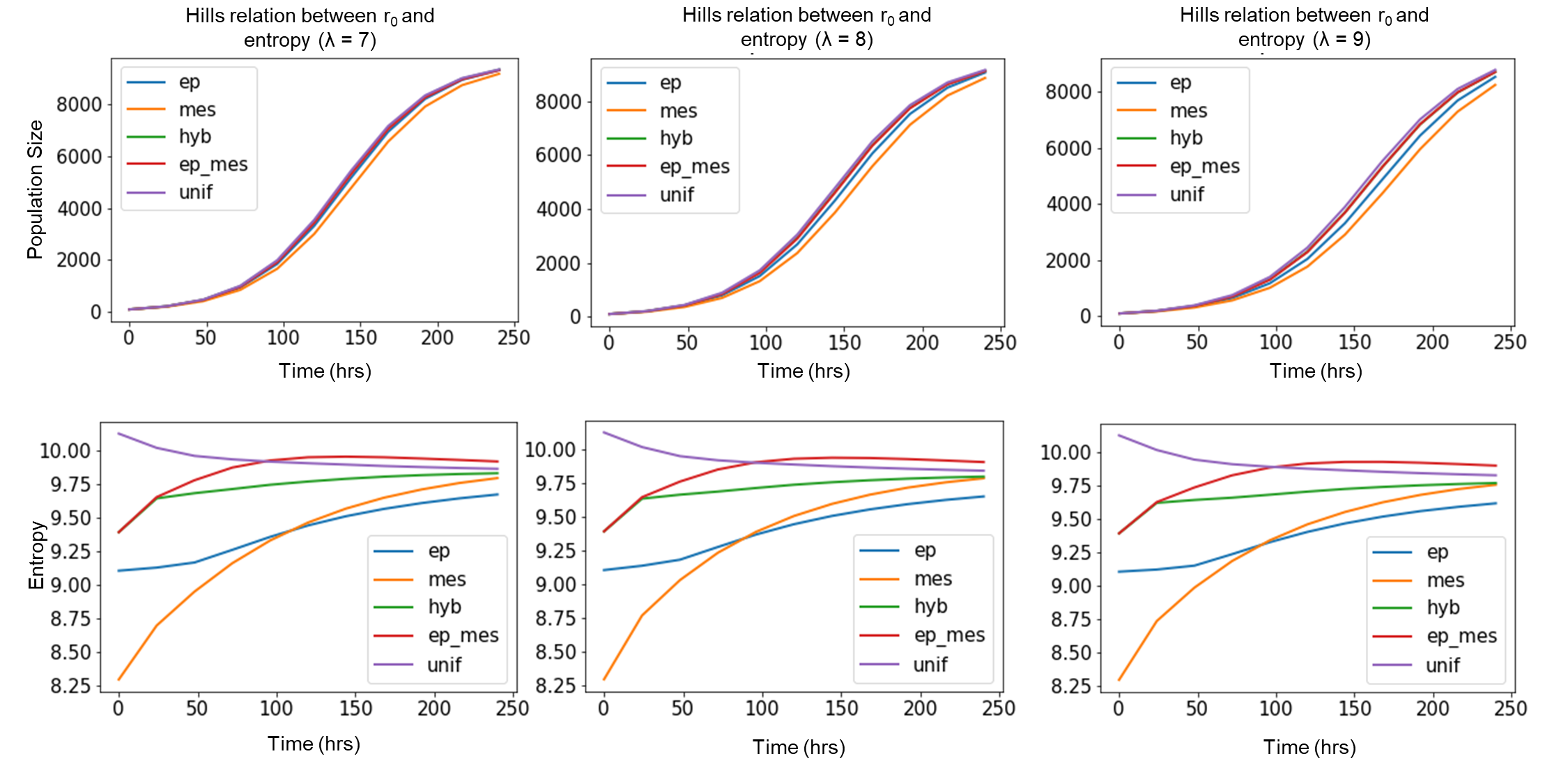}
    \caption{Temporal changes in total population size and population heterogeneity when growth rate of cell have sigmoidal relation with population heterogeneity (entropy). Parameters used to obtain above plot are $\alpha= 120$ hrs, $\eta_x = 1000$, and S levels in the population are set to 200K.}
    \label{figure S7}
\end{figure}

\section{Parameters for ODE \eqref{core regulation equation 2D} and \eqref{core regulation equation epigenetic regulation}}
\label{Annexes Parameters ODE}
We detail the parameters underlying ODE \eqref{core regulation equation 2D}:
\begin{equation*}
\begin{cases}
\dot \mu_{200}&=g_{\mu_{200}}H_{Z,\mu_{200}}(Z)H_{S,\mu_{200}}(S)-g_{m_Z}H_{Z,m_Z}(Z)H_{S,m_Z}(S)Q(\mu_{200})-k_{\mu_{200}}\mu_{200}\\
\dot Z&=g_Zg_{m_Z}H_{Z,m_Z}(Z)H_{S,m_Z}(S)P(\mu_{200})-k_ZZ\\
\end{cases}.
\label{}
\end{equation*}

Functions $H_{Z,\mu_{200}}$, $H_{S,\mu_{200}}$, $H_{Z,m_Z}$ and $H_{Z,m_Z}$ are shifted Hill functions which write under the form
\begin{align*}
H(X)=\frac{1+\lambda\big(\frac{X}{X_0}\big)^n}{1+\big(\frac{X}{X_0}\big)^n}. 
\end{align*}

The associated parameters are given in Table \ref{Table Hill functions}. 
\begin{table}[H]
\begin{center}

 \begin{tabular}{ |c c| c c |c c| c c| c c|} 
 \hline
  & & & & & Molecules & & Molecules.Hour$^{-1}$ & & Hour$^{-1}$\\
  \hline
  $n$ & $6$ & & & $\mu_0$ & 10K & & & & \\
  \hline
  $n_{Z, \mu_{200}}$ & $3$ & $\lambda_{Z, \mu_{200}}$ & $0.1$ & $Z^0_{\mu_{200}}$ & $220K$ & $g_{\mu_{200}}$ & $2.1K$ & $k_{\mu_{200}}$ & $0.05$  \\
  \hline
  $n_{S, \mu_{200}}$ & $2$ & $\lambda_{S, \mu_{200}}$ & $0.1$ & $S^0_{\mu_{200}}$ & $180K$ & $g_{Z}$ & $0.1K$ & $k_Z$ & $0.1$  \\
  \hline
    $n_{Z, m_Z}$ & $2$ & $\lambda_{Z, m_Z}$ & $7.5$ & $Z^0_{m_Z}$ & $25K$ & $g_{m_Z}$ & $11$ & $k_{m_Z}$ & $0.5$  \\
  \hline
    $n_{S, m_Z}$ & $2$ & $\lambda_{S, m_Z}$ & $10$ & $S^0_{m_Z}$ & $180K$ &  & &  &   \\
  \hline
\end{tabular}

\end{center}

\caption{Parameters for the Hill functions.}
\label{Table Hill functions}

\end{table}

The functions $Y_\mu$, $Y_m$ and $L$ are defined by
\begin{align*}
Y_\mu(\mu):=\sum\limits_{i=1}^{n}{i\gamma_{\mu_i}\binom{n}{i}M_n^i(\mu)}, \quad
Y_m(\mu):=\sum\limits_{i=1}^{n}{\gamma_{m_i}\binom{n}{i}M_n^i(\mu)}, \quad
L(\mu):=\sum\limits_{i=0}^{n}{l_{i}\binom{n}{i}M_n^i(\mu)},
\end{align*}
where $M_n^i:=\frac{(\mu/\mu_0)^i}{(1+\mu/\mu_0)^n}$. Finally functions $P$ and $Q$ are given by
\begin{align*}
    P(\mu):=\frac{L(\mu)}{Y_{m}(\mu)+k_{m_Z}}, \qquad
    Q(\mu):=\frac{Y_\mu(\mu)}{Y_m(\mu)+k_{m_Z}}
\end{align*}

Parameters for these three functions are given in Table \ref{Table other functions}.

\begin{table}[H]
\begin{center}
 \begin{tabular}{ |c|c|c|c|c|c|c|c| } 
 \hline
 i & 0 & 1 & 2 & 3 & 4 & 5 & 6 \\
 \hline
 $l_i$ & 1.0 & 0.6 & 0.3 & 0.1 & 0.05 & 0.05 & 0.05\\
 \hline
 $\gamma_{mi}$ &\ & 0.04 & 0.2 & 1.0 & 1.0 & 1.0 & 1.0 \\
 \hline
$\gamma_{\mu i}$ & \ & 0.005 & 0.05 & 0.5 & 0.5 & 0.5 & 0.5 \\
\hline
\end{tabular}
\end{center}
\caption{Parameters for the functions $L$, $Y_{\mu}$ and $Y_m$}
\label{Table other functions}
\end{table}

The additional parameters of ODE \eqref{core regulation equation epigenetic regulation} are $\alpha=0.15$, and $\beta(t)=240.$ when $S$ is non-decreasing, and $\beta(t)=720.$ when $S$ is decreasing.

\section{Reduction of the advection function}
\label{annexes reduction}

Let us introduce the segments, $I_{ep}=[150000, 185270.541082]$, $I_{ep-mes}=[185270.541082, 193286.5731462]$, $I_{ep-hyb-mes}=[193286.573146, 208817.635271]$, $I_{hyb-mes}=[208817.635271, 224649.298597]$, and $I_{mes}=[224649.298597, 250000]$, 
which correspond respectively to the values of $S$ for which the model is monostable (with a unique equilibrium point which corresponds to the epithelial phenotype), bistable (with two equilibrium points which correspond to the epithelial and the mesenchymal phenotypes), tristable, bistable (with two equilibrium points which correspond to the hybrid and the mesenchymal phenotypes), and monostable (with a unique equilibrium point which corresponds to the mensechymal phenotype), as can be seen on the bifurcation diagram (Figure \ref{figure 2} A). 

The values of stable equilibrium points representing the three possible phenotypes slightly varyn depending on the value of $S$ as illustrated by the bifurcation diagram. We approximate each of them via a five-order polynomial interpolation, \textit{i.e.} a polynomial of the form
\[a_5S^5+a_4S^4+a_3S^3+a_2S^2+a_1S+a_0.\]
$P_{mes}(S), P_{hyb}(S)$ and $P_{ep}(S)$ respectively approximate the mesenchymal, hybrid and epithelial phenotypes (three solid lines on the bifurcations diagrams), while $P_{u1}(S)$ and $P_{u2}(S)$ correspond to the unstable equilibrium points (two dotted lines in the bifurcation diagram). The coefficients defining these five polynomials are given in Table~\ref{Table Polynomial coefficients}.

%P_{ep}(S)=np.array([-1.98068372e-20,  1.72778767e-14, -6.01655669e-09,  1.04558220e-03,
% -9.08181301e+01,  3.18506176e+06])
%x_stable_center_coef=np.array([-9.97960480e-19,  1.04447036e-12, -4.37203191e-07,  9.14932022e-02,
% -9.57246289e+03,  4.00596473e+08])
%x_stable_right_coef=np.array([-6.10906482e-21,  6.84633959e-15, -3.06591965e-09,  6.85909353e-04,
% -7.66847769e+01,  3.43040294e+06])
%x_unstable_short_coef=np.array([ 1.98171027e-17, -1.99617838e-11,  8.04269750e-06, -1.62015993e+00,
%  1.63180639e+05, -6.57391347e+09])
%x_unstable_long_coef=np.array([ 2.50429596e-19, -2.55194809e-13,  1.03983466e-07, -2.11771508e-02,
%  2.15571288e+03, -8.77460219e+07])

\begin{table}[H]
\begin{center}
 %\begin{tabular}{ |c|c|c|c|c|c|c| } 
 %\hline
 % & $a_5$ & $a_4$ & $a_3$ & $a_2$ & $a_1$ & $a_0$  \\
 %\hline
 %$P_{mes}$ & $-6.109064e-21$ & $6.846339e-15$ & %$-3.065919e-09$ & $6.859093e-04$ & $ -7.668477e+01$ & %$3.430402e+06$ \\
%\hline
%$P_{u1}$ & $2.504295e-19$ & $-2.551948e-13$ &  %$1.039834e-07$ &  $-2.117715e-02$ & 
%  $2.155712e+03$ & $-8.774602e+07$\\
% \hline
% $P_{hyb}$ &$-9.979604e-19$&  $1.044470e-12$ & %$-4.372031e-07$ &  $9.149320e-02$ &
% $-9.572462e+03$ &  $4.005964e+08$\\
% \hline
% $P_{u2}$ &  $1.981710e-17$ & $-1.996178e-11$ &  %$8.042697e-06$ & $-1.620159e+00$ &
%  $1.631806e+05$ & $-6.573913e+09$\\
%\hline
%$P_{ep}$ & $-1.980683e-20$& $1.727787e-14$ & %$-6.016556e-09$ &  $1.045582e-03$ &
% $-9.081813e+01$  &$3.1850617e+06$\\
% \hline
%\end{tabular}

\begin{tabular}{|c|c|c|c|c|c|c|}
\hline
&   $P_{mes}$  &  $P_{u1}$ & $P_{hyb}$ & $P_{u2}$ & $P_{ep}$ \\
\hline
$a_5$ &  $-6.109064e-21$ & $2.504295e-19$ & $-9.979604e-19$&  $1.981710e-17$ & $-1.980683e-20$ \\
\hline
$a_4$ & $6.846339e-15$ & $-2.551948e-13$ & $1.044470e-12$ & $-1.996178e-11$ & $1.727787e-14$ \\
\hline
$a_3$ & $-3.065919e-09$ & $1.039834e-07$ & $-4.372031e-07$ & $8.042697e-06$ &$-6.016556e-09$ \\
\hline 
$a_2$ & $6.859093e-04$ & $-2.117715e-02$& $9.149320e-02$ &  $-1.620159e+00$ &  $1.045582e-03$ \\
\hline
$a_1$ &  $ -7.668477e+01$ &  $2.155712e+03$ & $-9.572462e+03$ &$1.631806e+05$ & $-9.081813e+01$ \\
\hline
$a_0$ &$3.430402e+06$ & $-8.774602e+07$& $4.005964e+08$& 
$-6.573913e+09$ & $3.1850617e+06$\\
\hline
\end{tabular}
\end{center}
\caption{Coefficients of the polynomials $P_{mes}$, $P_{u1}$, $P_{hyb}$, $P_{u2}$ and $P_{ep}$. }
\label{Table Polynomial coefficients}
\end{table}
%$P_{ep}$, $P_{hyb}$, $P_{mes}$, $P_{u2}$, $P_{u2}$ are three five-order polynomials which writes under the form 
%. The values of the coefficients in each cases are detailed in Table \ref{Table Polynomial coefficients}.  

We can finally define our one-dimensional reduced function:

\begin{itemize}
    \item If $S\in I_{ep}$: $\tilde f_r(x,S)=-(x-P_{ep}(S))$. 
    \item If $S\in I_{ep-mes}$: 
    \begin{itemize}
        \item If $x\leq 0.5(P_{mes}(S)+P_{u1}(S))$: $\tilde f_r(x,S)= -(x-P_{mes}(S))$
        \item If $x\in [0.5(P_{mes}(S)+P_{u1}(S)), 0.5(P_{u1}(S)+P_{ep}(S))]$, $ \tilde f_r(x,S)= x-P_{u1}(S)$
        \item If $x \geq 0.5(P_{u1}(S)+P_{ep}(S))$, $\tilde f_r(x,S)= -(x-P_{ep}(S))$
    \end{itemize}
    \item If $S\in I_{ep-hyb-mes}$:
    \begin{itemize}
        \item If $x\leq 0.5(P_{mes}(S)+P_{u1}(S))$, $\tilde f_r(x,S)= -(x-P_{mes}(S))$
        \item If $x\in [0.5(P_{mes}(S)+P_{u1}(S)), 0.5(P_{u1}(S)+ P_{hyb}(S))]$, $\tilde f_r(x,S)=(x-P_{u1}(S))$
        \item If $x\in [0.5(P_{u1}(S)+ P_{hyb}(S)), 0.5(P_{hyb}(S)+P_{u2}(S))]$, $\tilde f_r(x,S)=-(x-P_{hyb}(S))$
        \item If $x\in [0.5(P_{hyb}(S)+P_{u2}(S)), 0.5(P_{u2}(S)+P_{ep}(S))]$, $\tilde f_r(x,S)=x-P_{u2}(S)$
        \item If $x \geq 0.5(P_{u2}(S)+P_{ep}(S))$, $\tilde f_r(x,S)= -(x-P_{ep}(S))$
    \end{itemize}
    \item If $S\in I_{hyb-mes}$:
    \begin{itemize}
    \item If  $x\leq 0.5(P_{mes}(S)+P_{u1}(S))$, $f_r(x,S)= -(x-P_{mes}(S))$
    \item If $x\in [ 0.5(P_{mes}(S)+P_{u1}(S)), 0.5(P_{u1}(S)+P_{hyb(S)})]$, $\tilde f_r(x, S)=x-P_{u1}(S)$
    \item If $x\geq$, $\tilde f_r(x,S)= -(x-P_{hyb}(S))$
    \end{itemize}
    \item If $S\in I_{mes}$: $\tilde f_r(x,S)=-(x-P_{mes}(S))$
\end{itemize}

We are looking for a function that can be written as a multiple of $\tilde f_r$, \textit{i.e.} $f_r=k \tilde f_r$, with $k>0$. 
In order to choose the most suitable parameter, we look for the value of $k$ that minimises the quantity 
\begin{equation}
    \int_{0}^T\int_{150 K}^{250 K}\int_{0}^{25K}\int_0^{800K}{\lvert x(t, S, x_0)-\mu(t,S, x_0, Z_0) \rvert dt\,dS\,dx_0\, dZ_0},
    \label{minimised quantity fr}
\end{equation}
where for all $x_0 \in [0, 25K]$, $S\in [150K, 250K]$, $x(\cdot, S, x_0)$ solves
\begin{align*}
\begin{cases}
\dot x(t, S, x_0)=f_r\left(x(t, S, x_0)\right)\\
x(0,S,x_0)=x_0
\end{cases},
\end{align*}
and for all $Z_0\in [0, 800K]$, $\left(\mu(\cdot, S, x_0, Z_0), Z(\cdot, S, x_0, Z_0)\right)$ solves~\eqref{core regulation equation 2D}. 

In practice, this integral has been approximated for $T\in \{10, 100, 1000\}$ by the Riemann sum 
\[\frac{1}{N_T N_S N_x N_Z} \sum\limits_{i,j,k,l}{\lvert x(t_i, S_j, x0_{k})-\mu(t_i,S_j, x0_{k}, Z0_{l})\rvert },  \]
with $N_S=N_x=N_Z=20$, $N_T=100$, and for all $i,j,k,l$, $t_i=i \tfrac{T}{N_T}$, $S_j=150K+  j\tfrac{100K}{N_S}$, $x0_k=k \tfrac{25K}{N_x}$, and $Z0_l=l \tfrac{800K}{N_Z}$, where 
$x(\cdot, S_j, x0_k)$ solves
\begin{align*}
\begin{cases}
\dot x(t, S_j, x0_k)=f_r\left(x(t, S_j, x0_k)\right)\\
x(0,S_j,x0_k)=x0_k
\end{cases},
\end{align*}
$\left(\mu(\cdot, S_j, x0_k, Z0_l), Z(\cdot, S_j, x0_k, Z0_l)\right)$ solves~\eqref{core regulation equation 2D}. 

%In practice, this integral has been approximated, for $T\in \{10, 100, 1000\}$ with the sum {\color{red} add}

The value of the multiplicative constant which minimises~\eqref{minimised quantity fr} depends on $T$ but we establish that, for $T\in \{10, 100, 1000\}$ its value is rather insensitive to that of $T$: it is close to $0.02$, which it the value that we select. The obtained function is shown in Figure \ref{figure S2} for various values of $S$.

\bibliographystyle{vancouver}
\bibliography{biblio}

\begin{thebibliography}{10}

\bibitem{jacquemin2022Dynamic}
Jacquemin V, Antoine M, Dumont JE, Dom G, Detours V, Maenhaut C.
\newblock Dynamic Cancer Cell Heterogeneity: Diagnostic and Therapeutic
  Implications.
\newblock Cancers. 2022;14(2):280.
\newblock Available from: \url{https://doi.org/10.3390/CANCERS14020280}.

\bibitem{marusyk2020Intratumor}
Marusyk A, Janiszewska M, Polyak K.
\newblock Intratumor Heterogeneity: The Rosetta Stone of Therapy Resistance.
\newblock Cancer Cell. 2020;37(4):471–484.

\bibitem{bell2020Principles}
Bell CC, Gilan O.
\newblock Principles and mechanisms of non-genetic resistance in cancer.
\newblock British Journal of Cancer. 2020;122:465–472.

\bibitem{pillai2023Unraveling}
Pillai M, Hojel E, Jolly MK, Goyal Y.
\newblock Unraveling non-genetic heterogeneity in cancer with dynamical models
  and computational tools.
\newblock Nature Computational Science 2023 3:4. 2023;3(4):301–313.

\bibitem{brown2022Phenotypic}
Brown MS, Abdollahi B, Wilkins OM, Lu H, Chakraborty P, Ognjenovic NB, et~al.
\newblock Phenotypic heterogeneity driven by plasticity of the intermediate EMT
  state governs disease progression and metastasis in breast cancer.
\newblock Science Advances. 2022;8(31).

\bibitem{jain2022Population}
Jain P, Bhatia S, Thompson EW, Jolly MK.
\newblock Population Dynamics of Epithelial--Mesenchymal Heterogeneity in
  Cancer Cells.
\newblock Biomolecules. 2022;12(3):348.
\newblock Available from: \url{https://doi.org/10.3390/biom12030348}.

\bibitem{karacosta2019Mapping}
Karacosta LG, Anchang B, Ignatiadis N, Kimmey SC, Benson JA, Shrager JB, et~al.
\newblock Mapping Lung Cancer Epithelial-Mesenchymal Transition States and
  Trajectories with Single-Cell Resolution.
\newblock Nature Communications. 2019;10:5587.
\newblock Available from: \url{https://doi.org/10.1101/570341}.

\bibitem{sahoo2021mechanistic}
Sahoo S, Mishra A, Kaur H, Hari K, Muralidharan S, Mandal S, et~al.
\newblock A mechanistic model captures the emergence and implications of
  non-genetic heterogeneity and reversible drug resistance in ER+ breast cancer
  cells.
\newblock NAR Cancer. 2021;3(3).

\bibitem{font-clos2018Topography}
Font-Clos F, Zapperi S, Porta CAM.
\newblock Topography of epithelial–mesenchymal plasticity.
\newblock Proceedings of the National Academy of Sciences.
  2018;115(23):5902–5907.

\bibitem{hari2020Identifying}
Hari K, Sabuwala B, Subramani BV, Porta CAM, Zapperi S, Font-Clos F, et~al.
\newblock Identifying inhibitors of epithelial–mesenchymal plasticity using a
  network topology-based approach.
\newblock Npj Systems Biology and Applications 2020 6:1. 2020;6(1):1–12.

\bibitem{hong2015Ovol2}
Hong T, Watanabe K, Ta CH, Villarreal-Ponce A, Nie Q, Dai X.
\newblock An Ovol2-Zeb1 Mutual Inhibitory Circuit Governs Bidirectional and
  Multi-step Transition between Epithelial and Mesenchymal States.
\newblock PLOS Computational Biology. 2015;11(11):1004569.
\newblock Available from: \url{https://doi.org/10.1371/journal.pcbi.1004569}.

\bibitem{rashid2022Network}
Rashid M, Hari K, Thampi J, Santhosh NK, Jolly MK.
\newblock Network topology metrics explaining enrichment of hybrid
  epithelial/mesenchymal phenotypes in metastasis.
\newblock PLOS Computational Biology. 2022;18(11):1010687.
\newblock Available from: \url{https://doi.org/10.1371/JOURNAL.PCBI.1010687}.

\bibitem{steinway2015Combinatorial}
Steinway SN, Zañudo JGT, Michel PJ, Feith DJ, Loughran TP, Albert R.
\newblock Combinatorial interventions inhibit TGF$\beta$-driven
  epithelial-to-mesenchymal transition and support hybrid cellular phenotypes.
\newblock Npj Systems Biology and Applications. 2015;1:15014.
\newblock Available from: \url{https://doi.org/10.1038/npjsba.2015.14}.

\bibitem{george2017Survival}
George JT, Jolly MK, Xu S, Somarelli JA, Levine H.
\newblock Survival outcomes in cancer patients predicted by a partial EMT gene
  expression scoring metric.
\newblock Cancer Research. 2017;77(22):6415–6428.

\bibitem{ruscetti2016HDAC}
Ruscetti M, Dadashian EL, Guo W, Quach B, Mulholland DJ, Park JW, et~al.
\newblock HDAC inhibition impedes epithelial-mesenchymal plasticity and
  suppresses metastatic, castration-resistant prostate cancer.
\newblock Oncogene. 2016;35(29):3781–3795.

\bibitem{celia2018hysteresis}
Celi{\`a}-Terrassa T, Bastian C, Liu DD, Ell B, Aiello NM, Wei Y, et~al.
\newblock Hysteresis control of epithelial-mesenchymal transition dynamics
  conveys a distinct program with enhanced metastatic ability.
\newblock Nature communications. 2018;9(1):5005.

\bibitem{subbalakshmi2020nfatc}
Subbalakshmi AR, Kundnani D, Biswas K, Ghosh A, Hanash SM, Tripathi SC, et~al.
\newblock NFATc acts as a non-canonical phenotypic stability factor for a
  hybrid epithelial/mesenchymal phenotype.
\newblock Frontiers in oncology. 2020;10:553342.

\bibitem{hari2022Landscape}
Hari K, Ullanat V, Balasubramanian A, Gopalan A, Jolly MK.
\newblock Landscape of epithelial mesenchymal plasticity as an emergent
  property of coordinated teams in regulatory networks.
\newblock ELife. 2022;11.

\bibitem{boareto2016Notch}
Boareto M, Jolly MK, Goldman A, Pietil{\"a} M, Mani SA, Sengupta S, et~al.
\newblock Notch-Jagged signalling can give rise to clusters of cells exhibiting
  a hybrid epithelial/mesenchymal phenotype.
\newblock Journal of the Royal Society Interface. 2016;13(118):20151106.

\bibitem{jolly2017Inflammatory}
Jolly MK, Boareto M, Debeb BG, Aceto N, Farach-Carson MC, Woodward WA, et~al.
\newblock Inflammatory breast cancer: A model for investigating cluster-based
  dissemination.
\newblock Npj Breast Cancer. 2017;3(1):1–7.

\bibitem{neelakantan2017EMT}
Neelakantan D, Zhou H, Oliphant MUJ, Zhang X, Simon LM, Henke DM, et~al.
\newblock EMT cells increase breast cancer metastasis via paracrine GLI
  activation in neighbouring tumour cells.
\newblock Nature Communications. 2017;8:15773.
\newblock Available from: \url{https://doi.org/10.1038/ncomms15773}.

\bibitem{yamamoto2017Intratumoral}
Yamamoto M, Sakane K, Tominaga K, Gotoh N, Niwa T, Kikuchi Y, et~al.
\newblock Intratumoral bidirectional transitions between epithelial and
  mesenchymal cells in triple-negative breast cancer.
\newblock Cancer Science. 2017;108(6):1210–1222.

\bibitem{hitomi2021aAsymmetric}
Hitomi M, Chumakova AP, Silver DJ, Knudsen AM, Pontius WD, Murphy S, et~al.
\newblock Asymmetric cell division promotes therapeutic resistance in
  glioblastoma stem cells.
\newblock JCI Insight. 2021;6(3).

\bibitem{tripathi2020mechanism}
Tripathi S, Chakraborty P, Levine H, Jolly MK.
\newblock A mechanism for epithelial-mesenchymal heterogeneity in a population
  of cancer cells.
\newblock PLoS Computational Biology. 2020;16(2):1–27.

\bibitem{munsky2015Integrating}
Munsky B, Fox Z, Neuert G.
\newblock Integrating single-molecule experiments and discrete stochastic
  models to understand heterogeneous gene transcription dynamics.
\newblock Methods. 2015;85:12–21.

\bibitem{pally2022Extracellular}
Pally D, Goutham S, Bhat R.
\newblock Extracellular matrix as a driver for intratumoral heterogeneity.
\newblock Physical Biology. 2022;19(4):043001.
\newblock Available from: \url{https://doi.org/10.1088/1478-3975/AC6EB0}.

\bibitem{lovisa2015Epithelial}
Lovisa S, LeBleu VS, Tampe B, Sugimoto H, Vadnagara K, Carstens JL, et~al.
\newblock Epithelial-to-mesenchymal transition induces cell cycle arrest and
  parenchymal damage in renal fibrosis.
\newblock Nature Medicine. 2015;21(9):998–1009.

\bibitem{spencer2009non}
Spencer SL, Gaudet S, Albeck JG, Burke JM, Sorger PK.
\newblock Non-genetic origins of cell-to-cell variability in TRAIL-induced
  apoptosis.
\newblock Nature. 2009;459(7245):428--432.

\bibitem{strasen2018cell}
Strasen J, Sarma U, Jentsch M, Bohn S, Sheng C, Horbelt D, et~al.
\newblock Cell-specific responses to the cytokine TGF $\beta$ are determined by
  variability in protein levels.
\newblock Molecular systems biology. 2018;14(1):e7733.

\bibitem{jain2023epigenetic}
Jain P, Corbo S, Mohammad K, Sahoo S, Ranganathan S, George JT, et~al.
\newblock Epigenetic memory acquired during long-term EMT induction governs the
  recovery to the epithelial state.
\newblock Journal of the Royal Society Interface. 2023;20(198):20220627.

\bibitem{mantzaris2007single}
Mantzaris NV.
\newblock From single-cell genetic architecture to cell population dynamics:
  Quantitatively decomposing the effects of different population heterogeneity
  sources for a genetic network with positive feedback architecture.
\newblock Biophysical Journal. 2007;92(12):4271–4288.

\bibitem{hasenauer2011identification}
Hasenauer J, Waldherr S, Doszczak M, Radde N, Scheurich P, Allg{\"o}wer F.
\newblock Identification of models of heterogeneous cell populations from
  population snapshot data.
\newblock BMC bioinformatics. 2011;12:1--15.

\bibitem{shu2011Bistability}
Shu CC, Chatterjee A, Dunny G, Hu WS, Ramkrishna D.
\newblock Bistability versus bimodal distributions in gene regulatory processes
  from population balance.
\newblock PLoS Computational Biology. 2011;7(8).

\bibitem{spetsieris2009novel}
Spetsieris K, Zygourakis K, Mantzaris NV.
\newblock A novel assay based on fluorescence microscopy and image processing
  for determining phenotypic distributions of rod-shaped bacteria.
\newblock Biotechnology and Bioengineering. 2009;102(2):598–615.

\bibitem{hasenauer2012analysis}
Hasenauer J, Schittler D, Allg{\"o}wer F.
\newblock Analysis and simulation of division-and label-structured population
  models: a new tool to analyze proliferation assays.
\newblock Bulletin of mathematical biology. 2012;74:2692--2732.

\bibitem{schittler2013new}
Schittler D, Allg{\"o}wer F, De~Boer RJ.
\newblock A new model to simulate and analyze proliferating cell populations in
  BrdU labeling experiments.
\newblock BMC systems biology. 2013;7:1--6.

\bibitem{bhatia2019Interrogation}
Bhatia S, Monkman J, Blick T, Pinto C, Waltham A, Nagaraj SH, et~al.
\newblock Interrogation of phenotypic plasticity between epithelial and
  mesenchymal states in breast cancer.
\newblock J Clin Med. 2019;8(6):893.

\bibitem{chertock2017practical}
Chertock A.
\newblock A practical guide to deterministic particle methods.
\newblock In: Handbook of numerical analysis. vol.~18. Elsevier; 2017. p.
  177--202.

\bibitem{alvarez2023particle}
Alvarez FE, Guilberteau J.
\newblock A particle method for non-local advection-selection-mutation
  equations.
\newblock arXiv preprint arXiv:230414210. 2023.

\bibitem{lu2013microrna}
Lu M, Jolly MK, Levine H, Onuchic JN, Ben-Jacob E.
\newblock MicroRNA-based regulation of epithelial--hybrid--mesenchymal fate
  determination.
\newblock Proceedings of the National Academy of Sciences.
  2013;110(45):18144--18149.

\bibitem{jia2019possible}
Jia W, Deshmukh A, Mani SA, Jolly MK, Levine H.
\newblock A possible role for epigenetic feedback regulation in the dynamics of
  the epithelial-mesenchymal transition (EMT.
\newblock Physical Biology. 2019;16(6):066004.
\newblock Available from: \url{https://doi.org/10.1088/1478-3975/ab34df}.

\bibitem{sigal2006Variability}
Sigal A, Milo R, Cohen A, Geva-Zatorsky N, Klein Y, Liron Y, et~al.
\newblock Variability and memory of protein levels in human cells.
\newblock Nature. 2006;444(7119):643–646.

\bibitem{vega2004snail}
Vega S, Morales AV, Oca{\~n}a OH, Vald{\'e}s F, Fabregat I, Nieto MA.
\newblock Snail blocks the cell cycle and confers resistance to cell death.
\newblock Genes \& development. 2004;18(10):1131--1143.

\bibitem{corre2014stochastic}
Corre G, Stockholm D, Arnaud O, Kaneko G, Vi{\~n}uelas J, Yamagata Y, et~al.
\newblock Stochastic fluctuations and distributed control of gene expression
  impact cellular memory.
\newblock PLoS One. 2014;9(12):e115574.

\bibitem{nordick2022nonmodular}
Nordick B, Yu PY, Liao G, Hong T.
\newblock Nonmodular oscillator and switch based on RNA decay drive
  regeneration of multimodal gene expression.
\newblock Nucleic Acids Research. 2022;50(7):3693--3708.

\bibitem{loos2018hierarchical}
Loos C, Moeller K, Fr{\"o}hlich F, Hucho T, Hasenauer J.
\newblock A hierarchical, data-driven approach to modeling single-cell
  populations predicts latent causes of cell-to-cell variability.
\newblock Cell Systems. 2018;6(5):593--603.

\bibitem{mantzaris2006stochastic}
Mantzaris NV.
\newblock Stochastic and deterministic simulations of heterogeneous cell
  population dynamics.
\newblock Journal of theoretical biology. 2006;241(3):690--706.

\bibitem{degond1989weighted}
Degond P, Mas-Gallic S.
\newblock The weighted particle method for convection-diffusion equations. I.
  The case of an isotropic viscosity.
\newblock Mathematics of computation. 1989;53(188):485--507.

\bibitem{waldherr2018estimation}
Waldherr S.
\newblock Estimation methods for heterogeneous cell population models in
  systems biology.
\newblock Journal of The Royal Society Interface. 2018;15(147):20180530.

\bibitem{spetsieris2012single}
Spetsieris K, Zygourakis K.
\newblock Single-cell behavior and population heterogeneity: solving an inverse
  problem to compute the intrinsic physiological state functions.
\newblock Journal of biotechnology. 2012;158(3):80--90.

\bibitem{mantzaris2001numerical3}
Mantzaris NV, Daoutidis P, Srienc F.
\newblock Numerical solution of multi-variable cell population balance models.
  III. Finite element methods.
\newblock Computers \& Chemical Engineering. 2001;25(11-12):1463--1481.

\bibitem{mantzaris2001numerical1}
Mantzaris NV, Daoutidis P, Srienc F.
\newblock Numerical solution of multi-variable cell population balance models:
  I. Finite difference methods.
\newblock Computers \& Chemical Engineering. 2001;25(11-12):1411--1440.

\bibitem{dormand1980family}
Dormand JR, Prince PJ.
\newblock A family of embedded Runge-Kutta formulae.
\newblock Journal of computational and applied mathematics. 1980;6(1):19--26.

\end{thebibliography}

\end{document}